\documentclass[12pt]{article}
\usepackage{amsmath,amssymb,amsbsy,amsfonts,latexsym,amsopn,amstext,amsxtra,euscript,amscd,color}
\begin{document}
\bibliographystyle{plain}

\newfont{\teneufm}{eufm10}
\newfont{\seveneufm}{eufm7}
\newfont{\fiveeufm}{eufm5}
\newfam\eufmfam
              \textfont\eufmfam=\teneufm \scriptfont\eufmfam=\seveneufm
              \scriptscriptfont\eufmfam=\fiveeufm
%%%%%%%%%%%%%%%%%%%%%%%%%%%%%%%%%
%  \frak works on a single symbol at a time...
\def\frak#1{{\fam\eufmfam\relax#1}}

%%%%%%%%%%%%%%%%%%%  bbb-matter
\newcommand{\xor}[0]{\oplus}
\def\bbbr{{\rm I\!R}} %reelle Zahlen
\def\bbbm{{\rm I\!M}}
\def\bbbn{{\rm I\!N}} %natuerliche Zahlen
\def\bbbf{{\rm I\!F}}
\def\bbbh{{\rm I\!H}}
\def\bbbk{{\rm I\!K}}
\def\bbbp{{\rm I\!P}}
\def\bbbone{{\mathchoice {\rm 1\mskip-4mu l} {\rm 1\mskip-4mu l}
{\rm 1\mskip-4.5mu l} {\rm 1\mskip-5mu l}}}
\def\bbbc{{\mathchoice {\setbox0=\hbox{$\displaystyle\rm C$}\hbox{\hbox
to0pt{\kern0.4\wd0\vrule height0.9\ht0\hss}\box0}}
{\setbox0=\hbox{$\textstyle\rm C$}\hbox{\hbox
to0pt{\kern0.4\wd0\vrule height0.9\ht0\hss}\box0}}
{\setbox0=\hbox{$\scriptstyle\rm C$}\hbox{\hbox
to0pt{\kern0.4\wd0\vrule height0.9\ht0\hss}\box0}}
{\setbox0=\hbox{$\scriptscriptstyle\rm C$}\hbox{\hbox
to0pt{\kern0.4\wd0\vrule height0.9\ht0\hss}\box0}}}}
\def\bbbq{{\mathchoice {\setbox0=\hbox{$\displaystyle\rm
Q$}\hbox{\raise
0.15\ht0\hbox to0pt{\kern0.4\wd0\vrule height0.8\ht0\hss}\box0}}
{\setbox0=\hbox{$\textstyle\rm Q$}\hbox{\raise
0.15\ht0\hbox to0pt{\kern0.4\wd0\vrule height0.8\ht0\hss}\box0}}
{\setbox0=\hbox{$\scriptstyle\rm Q$}\hbox{\raise
0.15\ht0\hbox to0pt{\kern0.4\wd0\vrule height0.7\ht0\hss}\box0}}
{\setbox0=\hbox{$\scriptscriptstyle\rm Q$}\hbox{\raise
0.15\ht0\hbox to0pt{\kern0.4\wd0\vrule height0.7\ht0\hss}\box0}}}}
\def\bbbt{{\mathchoice {\setbox0=\hbox{$\displaystyle\rm
T$}\hbox{\hbox to0pt{\kern0.3\wd0\vrule height0.9\ht0\hss}\box0}}
{\setbox0=\hbox{$\textstyle\rm T$}\hbox{\hbox
to0pt{\kern0.3\wd0\vrule height0.9\ht0\hss}\box0}}
{\setbox0=\hbox{$\scriptstyle\rm T$}\hbox{\hbox
to0pt{\kern0.3\wd0\vrule height0.9\ht0\hss}\box0}}
{\setbox0=\hbox{$\scriptscriptstyle\rm T$}\hbox{\hbox
to0pt{\kern0.3\wd0\vrule height0.9\ht0\hss}\box0}}}}
\def\bbbs{{\mathchoice
{\setbox0=\hbox{$\displaystyle     \rm S$}\hbox{\raise0.5\ht0\hbox
to0pt{\kern0.35\wd0\vrule height0.45\ht0\hss}\hbox
to0pt{\kern0.55\wd0\vrule height0.5\ht0\hss}\box0}}
{\setbox0=\hbox{$\textstyle        \rm S$}\hbox{\raise0.5\ht0\hbox
to0pt{\kern0.35\wd0\vrule height0.45\ht0\hss}\hbox
to0pt{\kern0.55\wd0\vrule height0.5\ht0\hss}\box0}}
{\setbox0=\hbox{$\scriptstyle      \rm S$}\hbox{\raise0.5\ht0\hbox
to0pt{\kern0.35\wd0\vrule height0.45\ht0\hss}\raise0.05\ht0\hbox
to0pt{\kern0.5\wd0\vrule height0.45\ht0\hss}\box0}}
{\setbox0=\hbox{$\scriptscriptstyle\rm S$}\hbox{\raise0.5\ht0\hbox
to0pt{\kern0.4\wd0\vrule height0.45\ht0\hss}\raise0.05\ht0\hbox
to0pt{\kern0.55\wd0\vrule height0.45\ht0\hss}\box0}}}}
\def\bbbz{{\mathchoice {\hbox{$\sf\textstyle Z\kern-0.4em Z$}}
{\hbox{$\sf\textstyle Z\kern-0.4em Z$}}
{\hbox{$\sf\scriptstyle Z\kern-0.3em Z$}}
{\hbox{$\sf\scriptscriptstyle Z\kern-0.2em Z$}}}}
\def\ts{\thinspace}

\def\qed{\ifmmode
\squareforqed\else{\unskip\nobreak\hfil
\penalty50\hskip1em\null\nobreak\hfil\squareforqed
\parfillskip=0pt\finalhyphendemerits=0\endgraf}\fi}

\def\squareforqed{\hbox{\rlap{$\sqcap$}$\sqcup$}}

%%%%%%%%%%%%%%%%%%%%%%%%%
\def\cA{{\mathcal A}}
\def\cB{{\mathcal B}}
\def\cC{{\mathcal C}}
\def\cD{{\mathcal D}}
\def\cE{{\mathcal E}}
\def\cF{{\mathcal F}}
\def\cG{{\mathcal G}}
\def\cH{{\mathcal H}}
\def\cI{{\mathcal I}}
\def\cJ{{\mathcal J}}
\def\cK{{\mathcal K}}
\def\cL{{\mathcal L}}
\def\cM{{\mathcal M}}
\def\cN{{\mathcal N}}
\def\cO{{\mathcal O}}
\def\cP{{\mathcal P}}
\def\cQ{{\mathcal Q}}
\def\cR{{\mathcal R}}
\def\cS{{\mathcal S}}
\def\cT{{\mathcal T}}
\def\cU{{\mathcal U}}
\def\cV{{\mathcal V}}
\def\cW{{\mathcal W}}
\def\cX{{\mathcal X}}
\def\cY{{\mathcal Y}}
\def\cZ{{\mathcal Z}}
\newcommand{\rmod}[1]{\: \mbox{mod}\: #1}

\def\tcN{\cN^\mathbf{c}}
\def\Tr{{\mathrm{Tr}}}
\def\mand{\qquad \text{and} \qquad}
\renewcommand{\vec}[1]{\mathbf{#1}}
\def\udl#1{\underline{#1}}
\def\auth#1#2{\begin{color}{blue}{\it{#1,}}\end{color}\quad{\bf#2}:}
\def\Prob{{\item \quad }}
\def\Hint{{\bf Hint:}\quad}
\def\eqref#1{(\ref{#1})}

%%%%%%%%%%%%%%%%%%%%%%%%%%%%%%%%%%%
%%%%%%%%%%%%%%%%%%%%%%%%%%%%%%%%%%%

\newcommand{\ignore}[1]{}
\def\Dlog{\mathrm {Dlog}}
\hyphenation{re-pub-lished}
\parskip 1.5 mm
\def\lln{{\mathrm Lnln}}
\def\Res{\mathrm{Res}\,}

\def\F{{\bbbf}}
\def\Fp{\F_p}
\def\fp{\Fp^*}
\def\Fq{\F_q}
\def\ff{\F_2}
\def\ffn{\F_{2^n}}

\def\K{{\bbbk}}
\def \Z{{\bbbz}}
\def \N{{\bbbn}}
\def\Q{{\bbbq}}
\def \C{{\bbbc}}
\def \R{{\bbbr}}

\def\Zm{\Z_m}
\def \Um{{\mathcal U}_m}

\def \Bf{\frak B}

\def\Km{\cK_\mu}

\def\va {{\mathbf a}}
\def \vb {{\mathbf b}}
\def \vc {{\mathbf c}}
\def\vx{{\mathbf x}}
\def \vr {{\mathbf r}}
\def \vv {{\mathbf v}}
\def\vu{{\mathbf u}}
\def \vw{{\mathbf w}}
\def \vz {{\mathbfz}}

\def\({\left(}
\def\){\right)}
\def\fl#1{\left\lfloor#1\right\rfloor}
\def\rf#1{\left\lceil#1\right\rceil}

\def\flq#1{{\left\lfloor#1\right\rfloor}_q}
\def\flp#1{{\left\lfloor#1\right\rfloor}_p}
\def\flm#1{{\left\lfloor#1\right\rfloor}_m}

\def\Al{{\sl Alice}}
\def\Bob{{\sl Bob}}

\def\Or{{\mathcal O}}

\def\inv#1{\mbox{\rm{inv}}\,#1}
\def\invM#1{\mbox{\rm{inv}}_M\,#1}
\def\invp#1{\mbox{\rm{inv}}_p\,#1}

\def\Ln#1{\mbox{\rm{Ln}}\,#1}

\def \nd {\,|\hspace{-1.2mm}/\,}

\def\ord{\mu}

\def\E{\mathbf{E}}

\def\Cl{{\mathrm {Cl}}}

\def\epp{\mbox{\bf{e}}_{p-1}}
\def\ep{\mbox{\bf{e}}_p}
\def\eq{\mbox{\bf{e}}_q}

\def\bm{\bf{m}}

\newcommand{\floor}[1]{\lfloor {#1} \rfloor}

\newcommand{\comm}[1]{\marginpar{%
\vskip-\baselineskip %raise the marginpar a bit
\raggedright\footnotesize
\itshape\hrule\smallskip#1\par\smallskip\hrule}}

\def\rem{{\mathrm{\,rem\,}}}
\def\dist {{\mathrm{\,dist\,}}}

\def\veps{{\varepsilon}}
\def\eps{{\eta}}

\def\ind#1{{\mathrm {ind}}\,#1}
               \def \MSB{{\mathrm{MSB}}}
\newcommand{\abs}[1]{\left| #1 \right|}
%%%%%%%%%%%%%%%  Topmatter %%%%%%%%%%%%%%%%%%

\title{Divisibility, Smoothness
%%, Products of Small Primes And Their
and Cryptographic Applications}

\author{
 {\sc David Naccache}\\
 \'Equipe de cryptographie\\
  \'Ecole normale sup\'erieure\\
45 rue d'Ulm, F-75230 Paris, Cedex 05, France\\
 {\tt david.naccache@ens.fr}
\and
{\sc Igor E. Shparlinski}\\
             Department of Computing\\ Macquarie University\\
             Sydney,  NSW 2109, Australia\\
             {\tt igor@comp.mq.edu.au} }
\maketitle

\begin{abstract}
This paper deals with products of moderate-size primes, familiarly known as {\sl smooth numbers}. Smooth numbers play a crucial role in information theory, signal processing and cryptography.

We present various properties of smooth numbers relating to their enumeration, distribution and occurrence in various integer sequences. We then turn our attention to cryptographic applications in which smooth numbers play a pivotal role.
\end{abstract}

%% \paragraph{Keywords:}

\section{Introduction}

The goal of this paper is to shed light on the prominent role played by divisibility and smoothness in cryptography and related areas of
mathematics. This work intends to survey a wide range of results while steering away from too well-known examples. For doing so, we concentrate on some recently discovered applications of results about the arithmetic structure of integers.

We intend to convey to the reader a general comprehension of the state of the art, allow the devising of correct heuristics when problems cannot be tackled theoretically and help assessing the plausibility of new results.

In Section~\ref{sec:NT Background} we overview on a number of number-theoretic results commonly used for studying the multiplicative structure of integers. Most of the elementary results which we use are readily available from~\cite{HardyWright}; more advanced  results can be found, often in much more precise forms, in~\cite{CrPom,HalbRich,Harm2,IwKow,Ten1} and in many other standard analytic number theory manuals. Some of them are directly used in this paper, others remain in the background but we illustrate with them the variety of cryptographically useful analytic number theory tools.

We start our exploration of the worlds of divisibility and smoothness by asking a number of  natural questions. For instance, given a ``typical'' integer, what can be said about its largest divisor? Are Euler totient function values $\varphi(n)$ ``typical'' integers? What are the noteworthy properties of shifted primes $p-1$? How common are numbers who factor into products of primes which are all smaller than a bound $b$? The results listed here are neither exhaustive nor new (we refer the reader to references such number theory books or surveys such as~\cite{Granv10,HalTen,HilTen} for a more formal and systematic topical treatment).

Then, in Section~\ref{sec:Appl}, we use these results to shed light on a number of cryptographic constructions and attacks.

We remark that the specifics of this area is such that many impressive
works here may be underrated by
non-expetrs  as, at a first glance, they present only very small improvements over
previously known results. However these small improvements
are often principal steps forward and require the development of new ideas and very refined techniques.  Some  examples of such breakthrough
achievements include:

\begin{itemize}

\item  the estimate of Ford~\cite{Ford0}
on the counting function for the number of values of the Euler function,
see Section~\ref{sec:Euler};

\item the estimates of Ford~\cite{Ford2}
on the counting function of integers with an integer divisor in a
given interval, see Section~\ref{eq:div nat};

\item  the very tight estimates of
Croot, Granville, Pemantle \& Tetali~\cite{CGPT} on the stopping time
of the Dixon factoring algorithm, see Section~\ref{sec:Fact_DLog}.
\end{itemize}

Probably the oldest application of smoothness and divisibility is the celebrated {\sl Chinese Remainder Theorem} which allows 
%%igor us 
us to accelerate cryptographic functions and basic arithmetic operations using specific integer formats called {\sl residue number systems}, see~\cite{Mohan}. Results of this kind certainly deserve an independent treatment and we leave them outside of the framework of this paper.

We also recall that the idea of breaking a complex operation, depending on a parameter $n$, into a a recursion of simpler operations depending on the prime factors of $n$ can also be found in other fields such as signal processing.

The finite Fourier transform of a complex $n$-dimensional vector $y$ is the $n$-dimensional vector $Y$ defined as
$$
Y_k=\sum_{j=0}^{n-1} \omega_n^{jk}y_j.
$$
where  $\omega_n=e^{-2\pi i/n}$ is a complex $n$-th root of unity
Now, assume that $n=2m$ is even. We see that
$$
Y_k=\sum_{\substack{j=0 \\ j~\rm{even}}}^{n-1} \omega_n^{jk}y_j+\sum_{\substack{j=0 \\ j~\rm{odd}}}^{n-1}\omega_n^{jk}y_j
=\sum_{j=0}^{m-1} \omega_m^{jk}y_{2j}+\omega_n^{k}\sum_{j=0}^{m-1} \omega_m^{jk}y_{2j+1}
$$
In other words, the initial finite Fourier transform can be broken into two transforms of length $n/2$ on the projections of $y$ on even and odd dimensions. The same applies to divisibility by any prime or prime power and
%%igor us 
us  allows to derive a recursive Fast Fourier
Transform algorithm of sub-quadratic complexity when $n$ is smooth,~\cite{CoTu,vzGG}. Applications of this kind are also left out as we restrict
ourselves to the cryptographic genre.

%%igor

In general number theoretic results first appear as such, then, if appropriate, they are either directly applied or fine-tuned for cryptographic applications. Nonetheless, there are cases when important developments in cryptography have led to new arithmetic results.
For example, such is the bound of Coppersmith,  Howgrave-Graham \&
Nagaraj~\cite{CoH-GNag} on the number of divisors $d\mid n$ of a given integer $n$ in 
a prescribed arithmetic progression $d \equiv a \pmod k$, which is based on the ideas of the celebrated attack of Coppersmith~\cite{Copp1,Copp2} on RSA moduli with partially known factors. Yet another example is given by Boneh~\cite{Bon}, {sec:const smooth}
see Section~\ref{sec:const smooth} below.

\section{Conventions}

\subsection{Notations}
\label{sec:notation}

Throughout this paper we use Vinogradov's notation `$f(x) \ll g(x)$' which is equivalent to the Landau notation $f(x) = \cO(g(x))$, whilst being easier to chain as, for example, $f(x) \ll g(x) = h(x)$.\footnote{Note that $f(x) = \cO(g(x)) = h(x)$ is meaningless and $f(x) = \cO(g(x)) = \cO(h(x))$ may discard some useful information.} If convenient, we also write $g(x)\gg f(x)$ instead of $f(x)\ll g(x)$. We also write  $f(x)\asymp g(x)$ if $f(x)\ll g(x) \ll f(x)$.

The letter $p$ (possibly subscripted) always denotes a prime; $\varepsilon$ always stands for a small positive parameter on which implied constants
 may depend; $\log x$ denotes the natural logarithm of $x$. Calligraphic letters, for example, $\cA = \(a_n\)$, usually denote sequences of integers.

For a prime power $q$, we use $\F_q$ to denote the finite field of $q$ elements.

For an integer $m$, we use $\Z_N$ to denote the residue ring modulo $N$.

\subsection{Arithmetic Functions}

We use the following standard notations for the most common arithmetic functions for integers $m\ge 2$:

\begin{itemize}
\item $P(m)$, the largest prime divisor of $m$,
\item $\varphi(m)$, the Euler (totient) function of $m$,
\item $\omega(m)$, the number of distinct prime divisors of $m$,
\item$\tau(m)$,  the number of positive integer divisors of $m$.
\end{itemize}

Recall that $\varphi(m)$ is the number of positive integers $i \leq m$ with $\gcd(i,m)=1$
amdm that  $\tau(m)$ is sometimes denoted as $\sigma_0(m)$.

We also define $P(1) = \omega(1) = 0$ and $\tau(1) = \varphi(1) =1$.

Clearly $2^{\omega(m)}\leq \tau(m) $ and the inequality is tight for square-free $m \ge 1$.

Letting $x \ge 0$ be a real number, we denote by:
\begin{itemize}
\item $\pi(x)$     the number of primes $p\le x$,
%%igor \\
\item $\pi(x;q,a)$  the number of primes $p \le x$ such that  $p \equiv a \pmod q$.
\end{itemize}

\subsection{Integer Sequences}
\label{sec:Seq}

Besides the sequence of natural numbers $\N$, we devote in this paper particular attention to the following integer sequences:
\begin{itemize}
\item $\cP_a=\{ p+a \ :\  p\text{ prime}\}$,
\item $f(\N) = \{f(n) \ :\  n =1,2\ldots\}$,
\item $\varphi(\N) = \{ \varphi(n) \ :\  n =1,2\ldots\}$,
\item $\varphi(\cP_a) = \{ \varphi(p+a) \ :\  p\text{ prime}\}$.
\end{itemize}

In other words, $\cP_a$ is the sequence of shifted primes, $f(\N)$ is the sequence of polynomial valuations over $\N$, $\varphi(\N)$
is the sequence of Euler function values and $\varphi(\cP_a)$ is the sequence of Euler function values of shifted primes.

Amongst the sequence $\cP_a$, the instances $a =\pm 1$ are of special interest in cryptography and thus many papers concentrate only on these values. As results can usually be extended to any $a\ne 0$ at the cost of mere typographical changes, this work usually presents these results in this more general form.

\subsection{Smoothness}
\label{sec:def smooth}

$n\in\N$ is {\sl smooth} if $n$ has only {\sl small} prime divisors.
As the previous sentence does not define what {\sl small} is, we formally define $n$ as $y$-{\sl smooth} if all prime divisors $p \mid n$ are such that $p \le y$.

Alternatively, $n$ is $y$-smooth if and only if $P(n) \le y$.

Let $\cA = \(a_n\)$ be a sequence. We denote by $\psi(x,y;\cA)$ the number of $y$-smooth $a_n$ values found amongst the first $x$ elements of $\cA$ (that is, for $n \le x$). The following compact notations are used for the specific sequences defined in
Section~\ref{sec:Seq}:

\begin{eqnarray*}
\psi(x,y) &= & \psi(x,y;\N), \\
\pi_a(x,y)&= &  \psi(x,y;\cP_a),\\
\psi_f(x,y) &= &  \psi(x,y;f(\N)), \\
\Phi(x,y)&= &  \psi(x,y;\varphi(\N)),\\
\Pi_a(x,y)&= &  \psi(x,y;\varphi(\cP_a)).
\end{eqnarray*}

\subsection{The Dickman--de~Bruijn Function}

The {\sl Dickman--de Bruijn function} $\rho(u)$ is probably the  most popular smoothness density estimation tool.

$\rho(u)$ is defined recursively by:
$$\rho(u)=\left\{
\begin{array}{cl}
1 & \text{ if } 0 \leq u \leq 1,\\
\\
{\displaystyle 1 - \int_{n}^{t}\frac{\rho(v-1)}{v}dv} & \text{ if } u>1.\\
\end{array}
\right.
$$

Note that $\rho(u) = 1 - \log u$ for $1 \le u \le 2$.
For example, $\rho(\sqrt{e})=1/2$, that is, about half of the integers $n\le x$ has no prime divisors larger than $n^{1/\sqrt{e}} = n^{0.6065\ldots}$. This has been used by Vinogradov~\cite{Vin1} and by Burgess~\cite{Burg}, to estimate the smallest quadratic non-residue modulo a prime.

It is not difficult to show that as $u\to\infty$:
\begin{equation}
\label{eq:rho asymp}
\rho(u)=u^{-u +o(u)}
\end{equation}
and, more precisely,
$$
\rho(u)=\(\frac{e + o(1)}{u \log u}\)^u;
$$
even more accurate approximations to $\rho(u)$ are known,
see~\cite[Chapter~III.5, Theorem~8]{Ten1}.

\section{Number Theoretic Facts}
\label{sec:NT Background}

\subsection{Distribution of Primes}
\label{sec:PNT}

The {\sl Prime Number Theorem} states that for any fixed $A$:
\begin{equation}
\label{eq:PNT}
\pi(x) = {\mathrm {li}} x + \cO\(\frac{x}{(\log x)^A}\),
\end{equation}
where
$$
{\mathrm {li}}x = \int_{2}^x \frac{d\,t}{\log t}.
$$
Alternatively, using a more convenient (yet equivalent) formulation, in terms of the
$\vartheta$-function
$$
\vartheta(x) =\sum_{p \le x} \log p
$$
we can write that for any fixed $A>0$:
$$
\vartheta(x) = x + \cO\(\frac{x}{(\log x)^A} \).
$$

A commonly committed crime against primes is the assertion that:
$$
\pi(x) = \frac{x}{\log x} + \cO\(\frac{x}{(\log x)^A}\)
$$
 for any fixed $A>0$, which is {\sl wrong}, although, of course,
$$
 \pi(x) \sim {\mathrm {li}} x \sim \frac{x}{\log x}.
$$

An asymptotic estimate of the number of primes in arithmetic progressions is given by the {\sl Siegel--Walfisz theorem} see~\cite[Theorem~1.4.6]{CrPom} or, in an alternative form,~\cite[Chapter~II.8, Theorem~5]{Ten1}, which states that for every fixed  $A>0$ there exists $C>0$ such that for $x\ge 2$ and for all positive integers $q\le(\log x)^A$,
$$
\max_{\gcd(a,q)=1}\left|\pi(x;q,a)-\frac{{\mathrm {li}} x}{\varphi(q)}\right| \ll x \exp\(-C\sqrt{\log x}\),
$$
see also~\cite[Theorem~5.27]{IwKow}.

While for larger values of $q$, only conditional asymptotic formulae are known, for example, subject to the Generalized Riemann Hypothesis), the {\sl Brun-Titchmarsh theorem}, see~\cite[Chapter~3, Theorem~3.7]{HalbRich}, or~\cite[Theorem~6.6]{IwKow}, or~\cite[Chapter~I.4, Theorem~9]{Ten1}, gives a tight upper bound on $\pi(x;q,a)$ for all $q \le x^{1-\varepsilon}$. Namely, we have
$$
\pi(x;q,a) \ll \frac{x}{\varphi(q)\log (x/q)}
$$
without any restrictions on $x$ and $q$.

 Clearly for all $q \le x^{1-\varepsilon}$ we can replace  $\log(x/q)$ in the denominator with  $\log x$. Furthermore, this is conjectured to hold with just $\log x$ instead of $\log(x/q)$ in a wider range of $q$ (say up to $q \le x/(\log x)^A$ with some constant $A> 0$).

Finally, although for any given $q$, the Siegel--Walfisz theorem is the best know result, the {\sl Bombieri--Vinogradov theorem}, see~\cite[Theorem~17.1]{IwKow}, gives a much better estimate of $\pi(x;q,a)$ on {\sl average}  over $q$. In particular,
for every $A> 0$ there exists $B$ such that
$$
\sum_{q \le  \sqrt{x}/(\log x)^B}
\max_{y \le x} \max_{\gcd(a,q)=1}\left|\pi(y;q,a)-\frac{{\mathrm {li}} y}{\varphi(q)}\right|  \ll \frac{x}{(\log x)^A}.
$$

We conclude with the trivial but helpful remark that the bounds
$$
\pi(x;q,a) \le \pi(x) \mand
\pi(x;q,a) \le \frac{x}{q},
$$
can also be sufficient sometimes to establish useful results.

%% \subsection{Sums and Products of Primes}

\subsection{Mertens Formulae}

We recall the {\sl Mertens formulae} for the sums over primes
$$
\sum_{p\le x} \frac{1}{p}= \log \log x + A + o(1),  \qquad
\sum_{p\le x} \frac{\log p}{p} = \log x + B + o(1)
$$
and for the product
\begin{equation}
\label{eq:mert}
\prod_{p\le x} \(1 - \frac{1}{p}\) = \frac{C + o(1)}{\log x},
\end{equation}
where $A = 0.2614\ldots $, $B = 1.3325\ldots$, $C = e^\gamma = 1.7810\ldots$ and as before, $\gamma = 0.5772\ldots $ is the {\sl Euler-Mascheroni constant}, see~\cite[Sections~22.7 and~22.8]{HardyWright} or~\cite[Sections~I.1.4 and~I.1.5]{Ten1}. Vinogradov~\cite{Vinog-AI} gives a sharp bound of the error term.

Note that the  formula~\eqref{eq:mert} is related to the fact  that $\varphi(n)$ is rather large:
$$
n \ge \varphi(n) \gg  \frac{n}{\log \log n}.
$$

\subsection{Primes and the Zeta-Function}

The Riemann Zeta-function $\zeta(s)$ is defined for any $s\in\C$ with $\Re(s) > 1$ by
$$
\zeta(s) = \sum_{n=1}^\infty \frac{1}{n^s},
$$
and then  is analytically continued to all $s  \in\C$.

The {\sl Riemann Hypothesis} postulates that all the zeros of $\zeta(s)$ with $0 \le \Re(s) \le 1$ are such that $\Re(s)=1/2$. It is important to remind that there are other {\sl trivial} zeros outside of the {\sl critical strip} $0 \le \Re(s) \le 1$.

The {\sl Generalized Riemann Hypothesis} asserts that the same property holds   for a much wider class of similar functions called $L$-functions.

There are some {\sl explicit} formulae that relate $\pi(x)$ to the zeros of $\zeta$ in the critical strip.
In particular, the non-vanishing $\zeta(1 + it )\zeta(it) \ne 0$  for every $t \in \R$ implies the Prime Number Theorem under the form
$\pi(x) \sim {\mathrm {li}} x$. In fact the more we know about the distribution of the zeros of $\zeta$ the better
is the bound on $|\pi(x) - {\mathrm {li}} x|$ we get.

The best known result on the {\sl zero-free region} of $\zeta$ is due to Ford~\cite{Ford1}, who gives a more explicit version of the previous result obtained independently by Korobov~\cite{Kor} and Vinogradov~\cite{Vin2}, see also~\cite{IwKow}. In particular, thanks to these results, the
asymptotic formula~\eqref{eq:PNT} can be sharpened as
$$
\pi(x) - {\mathrm {li}} x  \ll   x \exp\(-C (\log x)^{3/5} ( \log\log x)^{-1/5}\)
$$
where $C = 0.2098$. A similar estimate for $\vartheta(x)$ can be obtained as well.

Unfortunately, besides the result of Ford~\cite{Ford1} and a few other similar estimates, very little progress has been witnessed in this area over the last decades.

For $\Re(s)> 1$, the {\sl Dirichlet product} is defined as:
\begin{eqnarray*}
\prod_{p}\(1-\frac{1}{p^s}\)^{-1} & = &
\prod_{p}\(1+\frac{1}{p^s}+ \frac{1}{p^{2s}}+ \frac{1}{p^{3s}}+  \ldots\)\\
& = & \sum_{n=1}^\infty \frac{1}{n^s} = \zeta(s).
\end{eqnarray*}

More generally, letting $\cS$ be any set of primes, and letting $\cN_\cS$ be the set of integers obtained by multiplying elements of $\cS$, we have:
\begin{equation}
\label{eq:Dir Prod}
\prod_{p \in \cS}\(1-\frac{1}{p^s}\)^{-1}= \sum_{n\in \cN_\cS} \frac{1}{n^s}.
\end{equation}

%%igor
\subsection{Beyond the Generalized Riemann Hypothesis}
\label{sec:Beyond GRH}

There is a common belief that the Generalized Riemann Hypothesis (GRH) fully characterizes the distribution of primes. This is unfortunately untrue and in many situations the GRH falls short of our expectations and heuristic predictions.
For example, for the gaps $d_n = p_{n+1}-p_n$ between consecutive primes $p_1 < p_2 < \ldots$ the GRH only implies that $d_n \ll p_n^{1/2}(\log p_n)^2$, while gaps are expected to be much smaller (and even be equal to $2$ infinitely often). Another example is the {\sl Elliott-Hallberstam Conjecture}, see~\cite[Section~17.1]{IwKow}, which asserts that for any fixed $\varepsilon > 0$ and $A>1$
$$
\sum_{q\le x^{1-\varepsilon}} \max_{\gcd(a,q)=1}
\left|\pi(x;q,a)-\frac{{\mathrm {li}} x}{\varphi(q)}\right|
\ll  \frac{x}{(\log x)^A}.
$$

On the other hand, and quite amazingly,  unconditional results on the distribution of primes which are stronger than results immediately implied by the GRH exist. One such estimates is the {\sl Brun-Titchmarsh theorem}, see Section~\ref{sec:PNT}. Other examples include a thread of works by Bombieri, Friedlander \& Iwaniec~\cite{BFI1,BFI2,BFI3} which extends Bombieri--Vinogradov's theorem, see Section~\ref{sec:PNT}, beyond the square-root range.

One of the
%%igor added "important" (for better formatting)
important
applications of these result is a remarkable result of Mikawa~\cite{Mik},
which asserts that for any fixed $a$ and almost all $q$
there is a prime $p \equiv a \pmod q$ with
%%igor display
$$
p \le q^{32/17 + o(1)}
$$
as $q\to \infty$. For all $q$, the best know estimate $p \ll q^{11/2}$ is due to Heath-Brown~\cite{H-B}.

\subsection{Euler Function}
\label{sec:Euler}

Here are a few beautiful properties of the Euler function which can be found in many standard number theory manuals (see, for example,~\cite{HardyWright}) For example, it is easy  to see that

\begin{equation}
\label{eq:phi div sum}
\sum_{d|m} \varphi(d)=m \mand \varphi(m) = m \sum_{d|m}\frac{\mu(d)}{d} ,
\end{equation}
where $\mu(d)$ is the M{\"o}bius function.
Furthermore, we have the identity
$$
\sum_{m=1}^{\infty}\frac{\varphi(m)q^m}{1-q^m}=\frac{q}{(1-q)^2}.
$$
Using~\eqref{eq:phi div sum} and simple analytic estimates
one can derive the following asymptotic formulae:
$$\frac{1}{m^2}\sum_{k=1}^m \varphi(k)=\frac{3}{\pi^2}+
\cO\(\frac{\log m}{m}\)
$$
and
$$
\frac{1}{m}\sum_{k=1}^m \frac{\varphi(k)}{k}=\frac{6}{\pi^2}+\cO\(\frac{\log m}{m}\).
$$

We also have explicit inequalities such as:
$$
\varphi(m)>\frac{m \log\log m}{e^\gamma (\log\log m)^2+3}
$$
for $m \ge 3$, where $\gamma = 0.5772\ldots $  is the {\sl Euler-Mascheroni constant}, and
$$
\varphi(m)>\sqrt{\frac{m}{2}}$$
for any $m \ge 1$.  Finally, for composite $m$:
$$
\varphi(m)\leq m-\sqrt{m}.
$$
There are also some much deeper questions about the Euler function.
One of the is studying the cardinality
$$
F(x) = \# \{\varphi(n) \le x\}
$$
of the set of values of the Euler function up to $x$,
for which Ford~\cite{Ford0} obtained a very precise estimate.

Ford~\cite{Ford0.5} has also established the validity of the
{\sl Serpi{\'n}ski conjecture} that for any integer $k\ge 2$
there is $m$ such that the equation $\varphi(n)=m$ has exactly $k$
solutions. We recall that by the {\sl Carmichael conjecture}
for any $m$ this equation has  either at least two solutions
or no solutions at all.

\section{How Smooth? How Many?}

\subsection{Empirical Estimates: A Cautionary Note}
\label{sec:Intuition}

Empirical estimates abound in cryptography. For examples, many cryptographers readily admit that, in the absence of obvious divisibility conditions, the density of primes in a given integer sequence is identical to the density of primes in $\N$. This and several similar ``postulates'' can be
frequently found throughout modern cryptographic literature. Let us illustrate the danger of such assumptions by a concrete example.

It is natural to approximate the probability that $p \nmid n$ when $n\le x$ is randomly chosen by $1 - 1/p$.

Now, assuming that all primes $p\le y$ are independent, we may infer that the probability that $p  \nmid  n$ for
all $x \ge p > y$ when $n\le x$ is chosen at random is close to:
$$
\prod_{x \ge p > y} \(1 - \frac{1}{p}\)=
\prod_{p \le x} \(1 - \frac{1}{p}\) \prod_{p \le y} \(1 - \frac{1}{p}\)^{-1}
\sim \frac{\log y}{\log x} = \frac{1}{u}
$$
by virtue of the Mertens formula, where $u$ is given by
\begin{equation}
\label{eq:Def u}
u = \frac{\log x}{\log y} \qquad \text{or}  \qquad  x = y^u.
\end{equation}

Here intuition leads to the seemingly elegant asymptotic formula
$$
\psi(x,y) \sim \frac{x}{u}.
$$
which is \ldots completely wrong!

\subsection{Estimating Smooth Integer Densities}

One of the most popular estimates of $\psi(x,y)$ is:
\begin{equation}
\label{eq:CEP bound}
\psi(x,y) = u^{-u + o(u)} x.
\end{equation}
This formula, due to Canfield, Erd{\H o}s \& Pomerance~\cite{CEP}, is applicable in the very large range:
$$u \le y^{1 - \varepsilon} \qquad \text{or} \qquad  y\ge (\log x)^{1+\varepsilon}$$
but the behavior of $\psi(x,y)$ changes for $y < \log x$.

While~\eqref{eq:CEP bound} is not an asymptotic formula (since $o(u)$ is in the exponent), asymptotic formulae for $\psi(x,y)$ exist. In particular, Hildebrand~\cite{Hild2} gave the asymptotic formula
\begin{equation}
\label{eq:Hild asymp}
\psi(x,y) \sim \rho(u) x
\end{equation}
for
$$u \le \exp \((\log y)^{3/5 - \varepsilon}\)  \qquad \text{or}
 \qquad
y \ge \exp\((\log \log x)^{5/3 + \varepsilon}\).
$$

A precise estimate of the error term in~\eqref{eq:Hild asymp} is given by Saias~\cite{Saias1}.

Note that~\eqref{eq:CEP bound} and~\eqref{eq:Hild asymp} imply~\eqref{eq:rho asymp}; of course~\eqref{eq:rho asymp} can also be obtained independently.

Unfortunately the  validity range of~\eqref{eq:Hild asymp} is much narrower than that of~\eqref{eq:CEP bound}, and is likely to  remain so for quite some time. Indeed, as per another result of Hildebrand~\cite{Hild1}, the validity of~\eqref{eq:Hild asymp} in the range:
$$
1\le u \le y^{1/2 - \varepsilon}  \qquad \text{or}
 \qquad
y\ge (\log x)^{2+\varepsilon}
$$
is equivalent to the Riemann Hypothesis.

\section{Estimating $\psi(x,y)$}

\subsection{Counting Very Smooth Numbers: Lattices}

To estimate $\psi(x,y)$ for very small values of $y$ one can resort to a geometric approach introduced by Ennola~\cite{Enn}, which has been
developed up to its natural limit by Granville~\cite{Granv1.5} (see also~\cite{Granv1}):

Let  $2 = p_1 < \ldots < p_s \le y$ be all $s = \pi(y)$ primes up to $y$. Then:
\begin{eqnarray*}
\psi(x,y) & = &\# \left\{(\alpha_1, \ldots, \alpha_s) \ :
\ \prod_{i =1}^s p_i^{\alpha_i} \le x\right\} \\
 & = &\# \left\{(\alpha_1, \ldots, \alpha_s) \ :
\ \sum_{i =1}^s \alpha_i \log p_i \le \log x\right\}.
\end{eqnarray*}

Thus our question boils-down to counting integer points in a specific tetrahedron. The number of integer points in any ``reasonable'' convex body is close to its volume. However, this is correct only if the volume is large with respect to its dimension $s$.

Thus we may expect that:
$$
\psi(x,y) \approx \frac{(\log x)^s }{s! \prod_{i =1}^s \log p_i}
$$
if $y$ is reasonably small. This approach can yield rigorous estimates,
see,  for example,~\cite{Granv1}.

\subsection{Upper Bounds: Rankin's Method}

For large values of $y$, the geometric approach fails to produce useful estimates. If only an upper bound is required, as is the case in many situations, then {\sl Rankin's method}~\cite{Ran} provides a reliable alternative.

Fix any constant $c > 0$.  Then
\begin{equation}
\label{eq:Rankin 1}
\psi(x,y) = \sum_{\substack{n \le x\\ p|n \Rightarrow p\le y}}
1 \le \sum_{\substack{n \le x\\ p|n \Rightarrow p\le y}}
\(\frac{x}{n}\)^c  \\
= \sum_{p|n \Rightarrow p\le y}
\(\frac{x}{n}\)^c.
\end{equation}

The underlying idea is that most of the contribution to $\psi(x,y)$ comes from integers which are close to $x$, so, although $(x/n)^c$ is larger than one for such integers, it is not much larger. On the other hand, $(x/n)^c$ decreases rapidly to zero when $n$ is much larger than $x$. So the above two steps do not cause over-counts.

 Using the fact that the right hand side of~\eqref{eq:Rankin 1} is an infinite series which can be represented as a Dirichlet product (see~\eqref{eq:Dir Prod}), we get:
\begin{equation}
\label{eq:Rankin 2}
\psi(x,y) \le  x^c   \sum_{p|n \Rightarrow p\le y}
 \frac{1}{n^c} =  x^c \prod_{p \le
y}\(1-\frac{1}{p^c}\)^{-1}.
\end{equation}

Using the Prime Number Theorem (in its best available asymptotic form) we estimate the product on the right hand side of~\eqref{eq:Rankin 2} as a function of $y$ and $c$ and minimize over all possible choices of $c>0$.

This task is technical but feasible and yields the quasi-optimal choice:
$$
c = 1 - \frac{u \log u }{\log y}
$$
which, in turn, yields an upper bound of the form~\eqref{eq:CEP bound}.

Simplicity (despite a few final technicalities) is the main advantage of this approach. In exchange, it suits only upper bounds and is apparently incapable of producing lower bounds.

\subsection{Asymptotic Formula: Buchstab--de Bruijn's  Recurrence}

We write each $y$-smooth $n>1$, as $n=pm$ where $p = P(n)$ is the largest prime factor of $n$. We note that $m \le x/p$ and is $p$-smooth.

Collecting together integers $n$ with $P(n) = p$ we get:
\begin{equation}
\label{eq:BdB}
\psi(x,y)= 1+\sum_{p \le y}\psi\(\frac{x}{p},p \)
\end{equation}
(where  $1$ at the front accounts for $n=1$), which is called the {\sl Buchstab--de Bruijn's  Recurrence}.

This recurrence formula has been used for both
lower and upper bounds and even for deriving asymptotic formulae.

We now use~\eqref{eq:BdB} to ``prove''~\eqref{eq:Hild asymp} for each fixed $u$
%%%igor
(we closely follow~\cite[Section~3.5]{Granv10}.

The ``proof'' is by induction over $N$, where $u \in (N, N+1]$. To ease the comprehension we deliberately ignore error terms and use the sign $\approx$ without specifying its formal meaning. However, we do guarantee to the reader that more careful analysis can re-cast the following formulae into a proper proof.

We start with the observation that for $0 < u \le 1$ we
trivially have $\psi(x,x^{1/u})  = \fl{x}$.

For $1 < u \le 2$ (that is, for $x \ge y \ge \sqrt{x}$), noticing that non-$y$-smooth numbers have one and only one prime divisor $p \ge y$, we get:

\begin{eqnarray*}
\psi(x,y) & = & x - \sum_{y \le p \le x } \#\{m : m \le \frac{x}{p}\}
=  x - \sum_{y \le p \le x } \fl{\frac{x}{p}}\\
& \approx & x - x \sum_{y \le p \le x } \frac{1}{p}
=  x\(1 -  \sum_{2\le p \le x } \frac{1}{p} + \sum_{2\le p \le y } \frac{1}{p}\).
\end{eqnarray*}
Therefore, by the Mertens formula,
\begin{eqnarray*}
\psi(x,y)
& \approx & x (1 - (\log \log x - \log \log y))  \\
& \approx & x \(1 -  \log \frac{\log x}{\log y}\)
= x (1 - \log u) = x \rho(u).
\end{eqnarray*}

We now note that the above step\ldots has not really necessary. It is nonetheless a good warming exercise for the next ``induction'' step.

Suppose that
$$\psi(x,x^{1/u}) \sim \frac{x}{\rho(u)}
$$
holds for $0 \le u \le N$.

Consider a  value of $u \in (N,N +1]$.

Subtracting the Buchstab--de Bruijn relation~\eqref{eq:BdB} with
$y=x^{1/N}$:
$$
\psi(x,x^{1/N})= 1+\sum_{p \le x^{1/N}}\psi\(\frac{x}{p},p \)
$$
from the same relation with $y=x^{1/u}$:
 $$
\psi(x,x^{1/u})= 1+\sum_{p \le x^{1/u}}\psi\(\frac{x}{p},p \).
$$
We obtain
\begin{eqnarray*}
\psi(x,x^{1/u})&=&\psi(x,x^{1/N})-
\sum_{x^{1/u}<p \le x^{1/N}}{\psi\(\frac{x}{p},p\)} \\ &\approx& x \(\rho(N) -
 \sum_{x^{1/u}< p \le x^{1/N}}
\frac{1}{p}\times\rho\(\frac{\log (x/p)}{\log p}\)\)
\end{eqnarray*}
since
$$
\frac{\log(x/p)}{\log p} = \frac{\log x}{\log p} - 1 < \frac{\log x}{\log (x^{1/u})} - 1 = u - 1 \le N,
$$
so the induction hypothesis applies (error terms ignored).

We now recall the definition of the function $\vartheta(z)$ and the Prime Number Theorem~\eqref{eq:PNT}. Writing $z = x^{1/t}$, by partial summation,  we get

\begin{eqnarray*}
 \sum_{x^{1/u}< p \le x^{1/N}}
\frac{1}{p}\rho\(\frac{\log (x/p)}{\log p}\)
&= &\int^{x^{1/N}}_{x^{1/u}} \rho \(\frac{\log x}{\log z} - 1 \)
 \frac{d\vartheta(z)}{z \log z} \\
&\approx& \int^{x^{1/N}}_{x^{1/u}} \rho \(\frac{\log x}{\log z} - 1 \)
\frac{dz}{z \log z}\\
&  = & \int^{u}_{N} \rho (t - 1) \frac{dt}{t}.
\end{eqnarray*}

We confess that a terrible offence has just been committed: instead differentiating $\vartheta(z)$ we have differentiated its approximation $z$. However, more careful examination shows that the above formulae are still correct.

Therefore
$$
\psi(x,x^{1/u}) \approx x\( \rho(N)  - \int^u_N \rho(t-1) \frac{dt}{t}\)
=   \rho(u)   x
$$
which concludes our ``proof''.

Estimating the largest prime divisor is a necessary step in many number-theoretic algorithms. For instance,  Bach, von zur Gathen \& Lenstra~\cite{BavzGaLe} introduce an algorithm for factoring polynomials over finite fields of characteristic $p$. The complexity of this algorithm depends on the largest prime divisor of the product the $k$ cyclotomic polynomial $\Phi_k(p)$ evaluated at $p$. A relationship between a number and its largest divisor allows to tune $k$ for every $p$ and optimize complexity. We refer the reader
to~\cite{Rony,Shoup,Zra2}
for more related results.

\section{Smoothness Miscellanea}

\subsection{Evaluating $\psi(x,y)$}

To optimize (balance) the complexity of steps in several cryptographic algorithms, one often needs more precise information about $\psi(x,y)$
than current\footnote{proven or conjectured} estimates and asymptotic formulae can provide.

For example, Parsell \& Sorenson~\cite{ParSor}, improving several previous results of Bernstein~\cite{Bern1}, have shown that for any parameter $\alpha$, one can estimate $\psi(x,y)$ up to a factor $1 + \cO\(\alpha^{-1} \log x\)$ in time
$$
\cO\(\frac{\alpha y^{2/3}}{\log y} + \alpha \log x\log \alpha\).
$$

A number of related results can be found in~\cite{HuntSor,Sor,Suz1,Suz2}.

\subsection{Constructing Constrained Smooth Numbers}
%%igor label added and referenced}
\label{sec:const smooth}

Producing smooth numbers is trivial.
However, constructing {\sl constrained} smooth numbers appears to be a challenging problem. A natural constraint, stemming from the study of digital signatures, is the requirement that the $y$-smooth number belongs to a given interval $[x, x+z]$.

Boneh~\cite{Bon}, motivated by certain cryptographic problems, has
devised a polynomial-time algorithm solving this problem for
some $x,y,z$ parameter combinations.

Results about the existence of very smooth numbers with a prescribed bit pattern at a certain position are given in~\cite{Shp1}, see also~\cite{GraShp} which gives an alternative approach (via character sums instead of exponential sums) that may probably be used to further improve the aforementioned result of~\cite{Shp1}.

More research in this area is certainly very desirable.

\subsection{Rough Numbers}

An integer $n$ is $y$-{\sl{rough}} if all prime divisors $p \mid n$ are such that $p > y$. We denote by $\Omega(x,y)$ be the number of $y$-{\sl{rough}} integers smaller than $x$.

Buchstab~\cite{Buch} gives the asymptotic formula:
$$
\Omega(x,y) \sim  \omega(u)\frac{x}{\log y} ,
$$
where the {\sl Buchstab} function $\omega(u)$ is defined as follows:
$$\omega(u)={\displaystyle \frac{1}{u}}\times\left\{
\begin{array}{cl}
1, & \quad \text{if } 1\leq u \leq 2,\\
\\
{\displaystyle 1 - \int_1^{u-1} \omega(t) dt}, &\quad \text{if } u\geq 2.\\
\end{array}
\right. $$

Rough numbers can be viewed  as ``approximations'' to primes. Rough numbers can be easily found and are proven to exist in various integer sequences of cryptographic interest. For example, Joye, Paillier \& Vaudenay~\cite{Joye1} use rough numbers as "interesting" candidates for primality testing during cryptographic key generation.

\subsection{Large Smooth Divisors}
\label{sec:Smooth Div}

It also natural to ask how often integers are expected to have a large smooth divisor; or, from a more quantitative perspective, explore the behavior of:
$$
\Theta(x,y,z) = \# \{ n \le x \ : \  \exists d \mid n, \ d > z, \
d\text{ is $y$-smooth}\}.
$$
While $\Theta$ has been addressed in the classical literature on smooth numbers, see~\cite{HalTen,Ten1,Ten2}, it has not  received as much attention as  $\psi(x,y)$.

Some asymptotic formulae for $\Theta$ have  recently been given by Banks \& Shparlinski~\cite{BaSh1} and Tenenbaum~\cite{Ten4}.
These formulae involve the parameter $u$ appearing in~\eqref{eq:Def u} and a parameter $v$ defined as:
$$
v =  \frac{\log z}{\log y}.
$$

The formulae also contain an integral involving $\rho(u)$ and its derivative.
Part of the motivation in~\cite{BaSh1} comes from a cryptographic problem discussed by Menezes~\cite{Men}, see as well Section~\ref{sec:HMQV}.

Shifted primes with large smooth divisors are studied by Pomerance \& Shparlinski~\cite{PomShp}.

\subsection{Next Largest Prime Divisors}

Characterizing the second largest prime divisor is of interest too, as the complexity of factoring an integer $n$ with Lenstra's elliptic curve factorization method~\cite{Len} (commonly called {\sl `the ECM'}), depends on this prime divisor.

More generally, denoting by $P_j(n)$ the $j$-th largest prime divisor of $n$, one may consider the joint distribution
$$
\psi(x,y_1, \ldots, y_k) =
\# \{ n \le x \ |\  P_j(n) \le y_j, \ j =1, \ldots, k\}.
$$

The work of Tenenbaum~\cite{Ten3} contains the most recent results and further references on this topic.

The case $k=2$ is especially important. Indeed,
using above notation, the ECM algorithm factors $n$ completely in time:
$$
\exp\(\(2+o(1)\) \sqrt{\log p \log \log p}\)n^{\cO(1)},
$$
where $p = P_2(n)$. This case has also got special attention
in~\cite{BachPer},
%%igor
see also~\cite{Suth} for some other applications.

\subsection{Other Facts}
\label{sec:mix}

In this section we present several unrelated results, which while unlikely to have any obvious cryptographic applications, still prove interesting for our exploration of smooth numbers.

Balog and Wooley~\cite{BalWool} have considered $k$-tuples of consecutive smooth integers and proved that for any $k$ and $\varepsilon > 0$ there are infinitely many $n$ such that $n+i$ is $n^\varepsilon$-smooth for $i = 1, \ldots, k$.
 In fact the proof in~\cite{BalWool} is based on very nice and
elementary explicit constructions.

One can also take  $k \to \infty$ and $\varepsilon \to 0$ (slowly) when  $n \to \infty$.

Balog~\cite{Bal2} proved that each sufficiently large integer $N$ can be written as $N= n_1 + n_2$ where $n_1, n_2$ are $N^{\alpha}$-smooth,
where
$$
\alpha = \frac{4}{9\sqrt e} = 0.2695\ldots\,.
$$

Results of this type may be considered as dual to the {\sl binary  Goldbach conjecture} claiming that all positive even integers $N \ge 4$ can be represented  as the sum of two primes.

Finally, various bounds of rational exponential sums
$$
S_{a,q}(x,y) = \sum_{\substack{n \le x \\n~\text{is $y$-smooth}}} \exp\(2 \pi  i \frac{a n}{q}\), \qquad \text{where }\gcd(a,q)=1,
$$
are given by Fouvry \& Tenenbaum~\cite{FouvTen1} and also by de~la~Bret{\`e}che \& Tenenbaum~\cite{BreTen}. Multiplicative character sums
$$
T_{a,q}(x,y) = \sum_{\substack{n \le x \\n~\text{is $y$-smooth}}} \chi(n-a), \qquad \text{where }\gcd(a,q)=1,
$$
with a nonprincipal multiplicative character $\chi$ modulo $q$ are estimated~\cite{Shp2}.
We also note that asymptotic formulae for the sums
$$
\sum_{\substack{a<n \le x\\
\\n~\text{is $y$-smooth}}}\frac{\varphi(n - a)}{n - a} \qquad\text{and}\qquad
 \frac{1}{\psi(x,y)}\sum_{\substack{a<n \le x\\
\\n~\text{is $y$-smooth}}} \varphi(n - a)
$$
are given in~\cite{LoiShp}. On the other hand, obtaining asymptotic
formulae (or even good estimates) for the sums
$$\frac{1}{\psi(x,y)} \sum_{\substack{a<n \le x\\
n~\text{is $y$-smooth}}}\tau(n - a) \qquad\text{and}\qquad
\frac{1}{\psi(x,y)}\sum_{\substack{a<n \le x\\
n~\text{is $y$-smooth}}}\omega(n - a),
$$
is still an open problem.

\section{Smoothness in Integers Sequences}

\subsection{Smooth Numbers in Arithmetic Progressions}

So far we considered the distribution of smooth values in the set of all natural numbers. A very natural generalization of this question, which is also of cryptographic interest, is the study of smooth numbers with an additional congruence condition.

In particular, we introduce the counting functions
$$
\psi(x,y;a,q) = \# \{ n \le x \ : \  n \text{ is $y$-smooth},\ n \equiv a \pmod q\}
$$
and
$$
\psi_q^*(x,y) = \# \{ n \le x \ : \  n \text{ is $y$-smooth},\ \gcd(n,q) = 1\}.
$$

Tenenbaum~\cite{Ten0} proved that
$$
 \psi_q(x,y) \sim \frac{\varphi(q)}{q} \psi(x,y)
$$
in a wide range of parameters.

In turn, a family of bounds of the forms
\begin{eqnarray*}
\psi(x,y;a,q) &\sim &\frac{1}{\varphi(q)} \psi_q(x,y), \\
\qquad
\psi(x,y;a,q) &\asymp & \frac{1}{\varphi(q)} \psi_q(x,y),\\
\qquad
\psi(x,y;a,q) &\gg &\frac{1}{\varphi(q)} \psi_q(x,y),
\end{eqnarray*}
(of decreasing strength but in increasingly larger ranges of $x$, $y$ and $q$) can be found in~\cite{BalPom,FouvTen2,Granv2,Granv3,Harm1,Sound}.
Some of these bounds hold for all $a$ with $\gcd(a,q)=1$, others hold only for almost all such integers $a$.

We note that bounds of exponential and character sums $S_{a,q}(x,y)$ and
$T_{a,q}(x,y)$, see Section~\ref{sec:mix}, can also be interpreted as  results about the uniformity of distribution of smooth numbers in arithmetic progressions ``on average''.

\subsection{Smooth Numbers in Small Intervals}
\label{sec:Smooth interv}

We now turn our attention to the sequence of integers in ``short'' intervals $[x, x+z)$.

Accordingly,
$$
\psi(x,y,z) = \psi(x+z,y) - \psi(x,y).
$$

It is natural to expect that:
$$
\psi(x,y,z) \sim \rho(u) z
$$
in a wide range of $x$, $y$ and $z$.

A series of ingenious  results due to Balog~\cite{Bal1},
 Croot~\cite{Croot2}, Friedlander \& Granville~\cite{FrGra}, Friedlander \& Lagarias~\cite{FrGra}, Harman~\cite{Harm1} and Xuan~\cite{Xuan} gives various interesting bits of information, but in general the status of
this problem is far from being satisfactory.

Croot's work~\cite{Croot2} is particularly interesting as it uses a quite unusual tool: bounds of {\sl bilinear Kloosterman sums} due to Duke, Friedlander \& Iwaniec~\cite{DFI}.

{From} both cryptographic and number theoretic perspectives,
the main challenge in that of obtaining good lower bounds in $\psi(x,y,4\sqrt{x})$, which appears to be currently out of reach.
This case is of special importance as it is crucial for the rigorous analysis of Lenstra's elliptic curve factoring algorithm~\cite{Len}. We note that the result
of Croot~\cite{Croot2} applies to intervals of similar length but unfortunately
for $y$ values which are much larger than these appearing in~\cite{Len}.

Finally, we recall that Lenstra, Pila \& Pomerance~\cite{LPP1,LPP2} have found an ingenious way to circumvent this problem by introducing a hyperelliptic factoring algorithm. For this algorithm, smooth numbers in large
intervals ought to be studied,  which
is already a feasible task. This has been achieved at the cost of very delicate arguments and required the developing of new algebraic and analytic tools by the authors.

\subsection{Smooth Shifted Primes}
\label{sec:smooth prime}

Let $\pi_a(x,y)$ be the counting function of smooth shifted primes
given in Section~\ref{sec:def smooth}.

It is strongly believed that for any fixed $a \ne 0$ the asymptotic formula
\begin{equation}
\label{eq:asymp shift primes}
\pi_a(x,y) \sim \rho(u) \pi(x)
\end{equation}
holds for a wide range of $x$ and $y$. Unfortunately, results of such breadth seem unreachable using current techniques.

However rather strong upper bounds are known. For example, Pomerance \& Shparlinski~\cite{PomShp} gave the estimate
$$
\pi_a(x,y)  \ll u \rho(u) \pi(x)
$$
for
$$
\exp\(\sqrt{\log x\log\log x}\)\le y \le x.
$$
In a shorter range
$$
\exp\((\log x)^{2/3+\varepsilon}\)\le y \le x
$$
the ``right'' upper bound
$$
\pi_a(x,y)  \ll \rho(u)\pi(x)
$$
follows from a result of Fouvry \& Tenenbaum~\cite[Theorem~4]{FouvTen2}.

It is not just the asymptotic formula~\eqref{eq:asymp shift primes} which is presently out of reach. In fact, even the obtaining of lower bounds
on $\pi_a(x,y)$ is an extremely difficult task where progress seems to be very slow.

The best known result, due to Baker \& Harman~\cite{BaHa} only asserts that there is a positive constant $A$ such that for $a \ne 0$,
$$
\pi_a(x,y)  \gg \frac{\pi(x)}{(\log x)^A}
$$
for $u \le 3.377\ldots$ (where as before, $u$ is defined by~\eqref{eq:Def u}), see also~\cite{Harm2}.

For most applications the logarithmic loss in the density of such primes is not important. However, if this becomes an issue, one can
use bound of Friedlander~\cite{Fried}:
$$
\pi_a(x,y) \gg \pi(x)
$$
which, however, is proven only for $u\le 2 \sqrt{e} = 3.2974\ldots$.

Finally, we recall yet another result of Baker \& Harman~\cite{BaHa} guarantees that
$$
\pi(x) - \pi_a(x,y) \gg  \pi(x)
$$
for $u \ge 1.477\ldots$.

The above results can be reformulated under the following equivalent forms which are usually better known and in which they are more frequently
used.

For some absolute constants $A,C > 0$ and for any $a\ne 0$:

\begin{itemize}
\item there are at least  $C \pi(x)/(\log x)^{A}$ primes $p\le x$ such that $p+a$ has a prime divisor $q \ge p^{0.6776}$;

\item there are at least  $C \pi(x)/(\log x)^{A}$  primes $p\le x$ such that all prime divisors $q$ of $p+a$ satisfy $q \le p^{0.2962}$.
\end{itemize}

The above two statements are expected to hold with $A = 0$, with $1 -\varepsilon$ instead of $0.6776$ and with $\varepsilon$ instead of $0.2962$ (for any $\varepsilon > 0$).

It is interesting to recall that results about shifted primes $p-1$ having a large prime divisor play a central role
in the deterministic primality test of Agrawal, Kayal \& Saxena~\cite{AgKaSa}.

\subsection{Smooth Values of Polynomials}

Let $f(X) \in \Z[X]$ and let $\psi_f (x,y)$ be defined as in Section~\ref{sec:def smooth}.

As in the case of shifted primes, rather strong upper bounds on $\psi_f(x,y)$ exist, see for example, the results of Hmyrova~\cite{Hm} and Timofeev~\cite{Tim}.

For a squarefree polynomial $f$, Martin~\cite{Mart} gives an asymptotic formula of the type
$$
\psi_f (x,y) \sim \rho(d_1u) \rho(d_2u) \dots \rho(d_k u) x,
$$
where $d_1, d_2, \ldots, d_k$ are the degrees of irreducible factors of $f$ over $\Z [x]$. This holds only  for very large values of $y$. See as well several related results by Dartyge, Martin \& Tenenbaum~\cite{DMT}, where smooth values of polynomials at prime valuations (that is, of $f(p)$) are also discussed.

\subsection{Smooth Totients}
\label{sec:smooth phi}

Let $\Phi(x,y)$ and $\Pi_a(x,y)$ be the functions counting smooth values of the Euler function on integers and shifted primes, respectively, see Section~\ref{sec:def smooth}.

Banks, Friedlander, Pomerance \& Shparlinski~\cite{BFPS} have shown
that in the range
$$
y \ge (\log\log x)^{1 + \varepsilon}
$$
we have
\begin{equation}
\label{eq:BFPS bound}
\Phi(x,y)\le x\exp(-(1+ o(1))\,u \log \log u).
\end{equation}

There are two interesting things to note about this result:
\begin{itemize}
\item the range is wider than that of~\eqref{eq:CEP bound}, with $\log \log y$ instead of $\log y$.

\item the bound is weaker than that of~\eqref{eq:CEP bound}, with $\log \log u$ instead of $\log u$ in the exponent.

\end{itemize}

Both bullet points reflect the fact that the values of Euler functions tend to be smoother than integers and shifted primes.
Furthermore, under a plausible conjecture about smooth shifted
primes similar to~\eqref{eq:asymp shift primes}, a
matching lower bound on $\Phi(x,y)$ has been obtained in~\cite{BFPS}.
Under even stronger conjectures, Lamzouri~\cite{Lamz} has obtained an
asymptotic formula for $\Phi(x,y)$ and similar quantities related to
iterations of the Euler function.

To estimate $\Pi_a(x,y)$, we can now use the trivial
inequality
$$
\Pi_a(x,y) \le \Phi(x+|a|,y),
$$ which, in fact, is quite sufficient for many applications. Furthermore, in a wide range of parameters, unless $u$ is small,
the above is equivalent to the expected estimate:
$$
\Pi_a(x,y) \ll   \pi(x) \exp\(-(1+ o(1))\,u \log \log u\)
$$
which is a full analogue of~\eqref{eq:BFPS bound}.

However for small values of $u$, obtaining the above estimate remains an important open problem.

On the other hand, in some ranges, several other bounds for $\Pi_a(x,y)$
have been obtained in~\cite{BFPS} using different techniques (such as a sieve method). For example, for
$y \ge \exp\(\sqrt{\log x\log\log x}\,\)$
we have
$$
\Pi_a(x,y) \ll \frac{\pi(x)}{u},
$$
and for $y \ge \log x $, we have
$$
\Pi_a(x,y) \le \frac{\pi(x)}{\exp\(\(1/2 + o(1)\) \sqrt{u} \log u\)}+\frac{\pi(x)\log\log x}{ \exp((1 + o(1)) u \log u) }.
$$

Note that estimates on the number of smooth values of the Euler function of polynomial sequences with integer and prime arguments, that is, $\varphi(f(n))$ and $\varphi(f(p))$, are given in~\cite{BFPS} as well.

\subsection{Smooth Cardinalities of Elliptic Curves}

One aspect of this question has already be mentioned in Section~\ref{sec:Smooth interv}
in relation to the elliptic curve factoring algorithm of Lenstra~\cite{Len}.

It is also interesting and important to study the arithmetic structure of
cardinalities of the reductions of a given elliptic curve defined over $Q$
modulo distinct primes. More precisely, given an elliptic curve  $\E$ over $Q$,
we denote by $N_p = \# \E(\F_p)$ the cardinality
of the set of rational points on the reduction of $\E$
modulo $p$ (for a sufficiently large prime $p$ such reduction always
leads to an elliptic curve over $\F_p$).

The number of prime divisors of $N_p$ has been studied
by Cojocaru~\cite{Coj2},
Iwaniec \&  Jim\'enez Urroz~\cite{IwJi-Ur},
Jim\'enez Urroz~\cite{Ji-Ur},
Liu~\cite{Liu1,Liu2,Liu3}, Miri \& Murty~\cite{MirMu} and
Steuding \& A. Weng~\cite{SteuWeng}.

Some heuristics about the number of prime values of $N_p$ for $p \le x$
has been discussed by Galbraith \& McKee~\cite{GalbMcK},
Koblitz~\cite{Kobl1,Kobl2} and  Weng~\cite{Weng}.
An upper bound on this quantity is obtained
by Cojocaru, Luca \& Shparlinski~\cite{CoLuSh},
see also~\cite{Coj1}.

However, it seems that there are no smootheness
results about the numbers $N_p$,
although this issue has been touched in McKee~\cite{McK}.
Probably obtaining an asymptotic formula or even a good lower
bound on the number of $y$-smooth values of $N_p$ for $p\le x$
is very hard, but perhaps some upper   bounds can be established.

Results of this kind are of great importance for the elliptic
curve cryptography.

\subsection{Smooth Class Numbers}

For a  integer $d<0$  we denote by $h(d)$ the
class number of the imaginary quadratic field $\Q(\sqrt{d})$.
Let $\cD$ be the set of {\sl fundamental discriminants},
that is,  the set of integers $d <0$ such that
\begin{itemize}
\item either $d \equiv 1 \pmod 4$ and $d$ is
square-free,
\item or  $d \equiv 0 \pmod 4$, $d/4 \equiv 2,3 \pmod 4$
and $d/4$ is square-free.
\end{itemize}

Using the so-called {\sl Cohen--Lenstra  heuristics} for
divisibility of class numbers,
see~\cite{CoLen}, Hamdy \& Saidak~\cite{HaSai} derived  a conditional asymptotic
formula for the number of
$d\in \cD$ with $-d \le x$ for which $h(d)$ is $y$-smooth.
Unfortunately, this seems to be
the only know result in this really exciting direction, see~\cite{BuchHa}
for the relevance of these problems to cryptography.
Studying the smoothness of  class numbers of other fields
is of great interest too
but is perhaps a very hard question.

\subsection{Smooth Numbers in Sumsets}

We recall that de~la~Bret{\`e}che~\cite{Bret1} gives a result of surprising generality and strength stating that, under some conditions, the proportion of smooth numbers among the sums $a+b$, where $a \in \cA$ and $b \in \cB$ is close to the expected value for a wide class of sets $\cA, \cB \subseteq \Z$.

For example, let $\cA$ and $\cB$ be two sets of integers in the interval $[1,x]$. Then, for any fixed $\varepsilon > 0$ and
uniformly for
$$
\exp\((\log x)^{2/3+\varepsilon}\)<y\le x,
$$
we have
\begin{equation*}
\begin{split}
&\#\{(a,b)\in\cA\times\cB~:~ a + b \ \text{is $y$-smooth}\}\\
&\qquad\qquad\qquad= \rho(u)\cdot\# \cA \# \cB
\(1 + \cO\(\frac{x\log(u + 1)}
{\sqrt{\# \cA \# \cB}\log y}\)\),
\end{split}
\end{equation*}
where, as usual, $u$ is given by~\eqref{eq:Def u}.

Although the authors are unaware of any immediate
cryptographic applications of
this result, we underline its high potential, given its generality and
``condition-free'' formulation.

Several other relevant results are given by Croot~\cite{Croot1}.

\subsection{Smooth Polynomials Over Finite Fields}
\label{sec:smooth poly}

In full analogy with the case of integers, we say that a polynomial $F \in \K[x]$ over a  field $\K$ is $k$-smooth if all irreducible
divisors $f \mid F$ satisfy $\deg f  \le k$.

For a finite field $\F_q$ of $q$-elements we denote
$$
N_q (m,k) = \# \{ f \in \F_q[x] \ : \  \deg f \le m \ f \text{ is $k$-smooth and monic}\}.
$$
Define
$$
u = \frac{m}{k} = \frac{\log q^m}{\log q^k}
$$
(the last expression makes the analogy with formula~\eqref{eq:Def u} completely explicit).

The systematic study of $N_q(m,k)$ dates back to the work of Odlyzko~\cite{Odl} who also discovered the relevance
of this quantity to the discrete logarithm problem in finite fields.

Several very precise results about $N_q(m,k)$ have recently been given by Bender \& Pomerance~\cite{BenPom}.
For example, by~\cite[Theorem~2.1]{BenPom} we have
$$
N_q(m,k) = u^{-u + o(u) }q^m
$$
as $k \to \infty$ and $u \to \infty$, uniformly for $q^k \ge m (\log m)^2$,
 and by~\cite[Theorem~2.2]{BenPom} we also have
$$
N_q(m,k)\ge \frac{q^m}{m^u}
$$
for $k \le \sqrt{m}$.

\section{Distribution of Divisors}

\subsection{More On Intuition}

We have already seen in Section~\ref{sec:Intuition} that carelessly applied intuition may lead to wrong conclusions. The following is yet another example.

It is obvious that the density of perfect squares $n = d^2$ is extremely small as there are only $\lfloor\sqrt{x}\rfloor$ perfect squares up to $x$.

Let us relax the relation $n = k^2$ and consider $n = km$ with  $k \le m \le k^{1.001}$ for $k,m\in N$. Such integers can be called quasi-squares.

It is natural to ask whether the density of quasi-squares is still small.
Say, are there only $o(x)$ quasi-squares up to $x$?
As the below results indicate, such quasi-squares form a set
of positive density which perhaps does not match the
intuition (too bad for the intuition \ldots).

\subsection{Notations}

Given a sequence of integers $\cA = \(a_n\)$, we denote
$$
H(x,y,z;\cA) = \#\{ a_n\le x \ :\  \exists \, d|a_n \text{ with } y<d\le
z\}.
$$
As usual, in the case of $\cA = \N$ we define
$$
H(x,y,z) = H(x,y,z;\N).
$$

\subsection{Natural Numbers}
\label{eq:div nat}

This case goes back to two old questions of Erd{\H o}s:
\begin{quotation}
{\sl Given an integer $N$ what is the size $M(N)$ of the {\sl multiplication table} $\{nm~:~ 1\le m,n \le \sqrt{N}\}$?}
\end{quotation}
and
\begin{quotation}
{\sl Is it true that almost all integers $n$ have
two divisors $d_1\mid n$ and $d_2\mid n$ with $d_1 < d_2 < 2d_1$?\/}
\end{quotation}

Hall \& Tenenbaum's book~\cite{HalTen} contains a very detailed treatment of such questions. For example, for the size of the multiplication table we have
from~\cite[Corollary~3]{Ford2}
$$
M(N) \asymp \frac{N}{ (\log N)^{\delta} (\log\log N)^{3/2}}
$$
as $N \to \infty$, where
\begin{equation}
\label{eq:Erd Number}
\delta = 1 - \frac{1 + \log \log 2}{\log 2} = 0.008607 \ldots
\end{equation}
is the {\sl Erd{\H o}s number}, see also~\cite[Theorem~23]{Ten1} for a slightly less precise result.

Ford~\cite{Ford2} has recently obtained a series of remarkable improvements of
several previously know results.
Several related results can also be found in~\cite{Bret2,Ford3,FordTen,Kou,Ten1}.
Unfortunately, exact formulations of precise results lead to rather cluttered technical conditions and estimates, which also depend on the relative sizes of $x$, $y$ and $z$ as well as on $z-y$ and $z/y$. Thus we limit our discussion to only a few sample results.

For example, let us define $v> 0$ by the relation
$$
z = y^{1 + 1/v}.
$$
Then, by~\cite[Theorem~1]{Ford2} for any real $x$, $y$ and $z$
with
$$
x \ge \max\{100000,y^2\} \mand
x \ge z \ge y \ge  100 $$
 we have
$$
\frac{H(x,y,z,\N)}{x} \asymp
\left\{
\begin{array}{cl}
v^{\delta} (\log v)^{-3/2} & \text{ if } 2y \le z \le y^{2},\\
\\
1 & \text{ if } z \ge y^2,
\end{array}
\right.
$$
where $\delta$ is the Erd{\H o}s number.

In particular, we see that for $\varepsilon>0$ and any sufficiently large $y$, we have
\begin{equation}
\label{eq:arb int}
H(x,y,y^{1 + \varepsilon},\N) \gg x,
\end{equation}
where the implied constant only depends on $\varepsilon>0$. Thus, there is a positive density of integers $n \le x$, depending only on $\varepsilon>0$, which have a divisor $d\mid n$ in the interval $d \in [y, y^{1+\varepsilon}]$.

We now prove~\eqref{eq:arb int} in the special case where $y$ is a power of $x$. That is, we prove that for $0 < \alpha < \beta < 1$:
\begin{equation}
\label{eq:large int}
 H(x,x^{\alpha},x^{\beta},\N) \gg x.
\end{equation}

In our proof we consider only {\sl prime divisors} $p \in [x^{\alpha},x^{\beta}]$ (instead of integer divisors) and
make the following two trivial observations:
\begin{itemize}
\item there are  $x/p + \cO(1)$
integers $n\le x$ divisible by $p$;
\item each $n \le x$ may have at most $K = \rf{1/\alpha}$ of them.
\end{itemize}
Hence,
$$
 H(x,x^{\alpha},x^{\beta},\N) \ge \frac{1}{K} \sum_{x^{\alpha} \le p \le x^{\beta}} \(\frac{x}{p} + \cO(1)\)
$$
and the sum on the right hand side
counts every integer $n \le x$ with a prime divisor $p \in [x^{\alpha},x^{\beta}]$ at most $K$ times. Therefore
%%igor remove "="
$$
 H(x,x^{\alpha},x^{\beta},\N) \ge
%% =
\frac{x}{K} \sum_{x^{\alpha} \le p \le x^{\beta}} \frac{1}{p} + \cO(x^b).
 $$
By the Mertens formula, we now obtain
 \begin{eqnarray*}
 H(x,x^{\alpha},x^{\beta},\N) & \ge &
 \frac{x}{K} \(\log \log (x^{\beta}) -  \log \log (x^{\alpha}) + o(1)\)  \\
& = &  \frac{x}{K} \(\log \frac{\log (x^{\beta})}{\log (x^{\alpha})} + o(1)\) \\
& = &  \frac{1}{K} \(\log \left(\frac{\beta}{\alpha}\right) + o(1)\) x
\end{eqnarray*}
and~\eqref{eq:large int} follows.

There are other tell tale signs that integer divisors are densely distributed. For example, for an integer $s\ge 1$ we denote
$$
T(n) =  \max_{i=1,\ldots,\tau(n)-1} \frac{d_{i+1}}{d_i},
$$
where $1= d_1<\ldots<d_{\tau(n)}=n$ are the positive divisors of $n$.
Clearly,
$$
T(n) \le  P(n).
$$
However, for many integers $T(n)$ is much smaller than $P(n)$. By a result of Saias~\cite[Theorem~1]{Saias2}, we know that
for any fixed $t$ and sufficiently large $x$,
$$
\# \left\{n \le x\ : \  T(n) \le t \right\} \asymp
\frac{x\log t}{\log x}.
$$

\subsection{Shifted Primes}
\label{sec: div shift p}

Ford~\cite{Ford2} has given upper bounds on $H(x,y,z;\cP_a)$ of about the same strength as these applying to $H(x,y,z;\cN)$, where $\cP_a$ is
defined in Section~\ref{sec:Seq}.

The situation with lower bounds $H(x,y,z;\cP_a)$ is quite bleak, although heuristically there is little doubt that $H(x,y,z;\cP_a)$
should behave similarly to $H(x,y,z;\cN)$.

One of the very few known lower bounds (yet, with many cryptographic applications) is given in~\cite[Theorem~7]{Ford2}:
for $a \ne 0$ and $0 < \alpha < \beta$:
$$
 H(x,x^{\alpha},x^{\beta},\cP_a) \gg  \pi(x)
$$
(where the implied constant depends on $a$, $\alpha$ and $\beta$).

The proof is similar to our proof of~\eqref{eq:large int}, but requires some technical analytic number theory tools, namely, the Bombieri--Vinogradov  theorem,
see Section~\ref{sec:PNT},  since instead of {\sl integers} $n \le x$ with $p\mid n$ we need to count {\sl primes} $q \le x$ with $p \mid (q-a)$.

To implement this approach one also needs the elementary observation that is enough to consider only the case $0 < \alpha <  \beta \le 1/2$
%%igor $d|n$
(since if $d\div n$ then
%%igor $(n/d) \mid n$).
$n/d$ is also a divisor of $n$).

%%igor
Finally, we remark that the Brun pure sieve (that is, a properly truncated version of the inclusion-exclusion principle), see~\cite[Theorem~2.3]{HalbRich}, or~\cite[Theorem~3, Section~I.4.2]{Ten1}, immediately implies that for any $a \ne 0$
\begin{equation}
\label{eq:Div p-1-Brun}
H(x,y,z;\cP_a) = \(1 + O\(\frac{\log y}{\log z}\)\) \pi(x).
\end{equation}

For example, the bound~\eqref{eq:Div p-1-Brun} can be used for the analysis of some cryptographic attacks Cheon~\cite{Cheon}.

\subsection{Polynomials}

Unfortunately there seem to be no results about the distribution of integers divisors of polynomials. Nonetheless, this problem does not look hopeless.

%%igor
\subsection{Cardinalities of Elliptic Curves}

As in the case of polynomials, there seem to be no results about the distribution of divisors of cardinalities of elliptic curves over finite fields. The question is certainly hard but not completely hopeless and deserves to be studied. Furthermore, Menezes \& Ustooglu~\cite{MenUst} point out that
this question has direct cryptographic applications.

It is probable that for the set of all elliptic curves over a given finite field $\F_q$ new results can be obtained by combining the Brun sieve technique, see~\cite[Theorem~2.3]{HalbRich}, or~\cite[Theorem~3, Section~I.4.2]{Ten1} with results of Howe~\cite{Howe} on divisibility statistics of elliptic curves.

\subsection{Totients}

Here is another confirmation that totients are not typical integers.

As we have mentioned, $H(x,y,z;\cP_a)$ is expected to behave similarly to $H(x,y,z;\cN)$. However the behaviour of $H(x,y,z;\varphi(\N))$ is very different.

Given that typical values of the Euler function
\begin{itemize}
\item have more prime divisors, due to a result of Erd{\H o}s \& Pomerance~\cite{ErdPom},
\item have more integer divisors, due to a result of Luca \& Pomerance~\cite{LucPom},
\item are smoother, due to a result of Banks, Friedlander,
Pomerance \& Shparlinski~\cite{BFPS},
see also Section~\ref{sec:smooth phi},
\end{itemize}
than a typical integer, it is also natural to expect that totients have denser divisor sets. This is supported by several recent
results of Ford \& Hu~\cite{FordHu}, who in particular show that

\begin{itemize}
\item uniformly over $1 \le y \le x/2$, we have $H(x,y,2y;\varphi(\N)) \gg x$;
\item for  $y = x^{o(1)}$, we have $H(x,y,2y;\varphi(\N)) \sim x$;
\item for a positive proportion of integers $n$, there is a divisor $d \mid \varphi(n)$ in every interval of the form $[K, 2K]$, $1 \le K \le n$.
\end{itemize}

\section{Cryptographic Applications}
\label{sec:Appl}

\subsection{Smoothness in Factoring and Discrete Logarithms}
\label{sec:Fact_DLog}

Most integer factorization algorithms, such as Dixon's method, the Quadratic Sieve, index calculus, the Number Field Sieve or Elliptic Curve Factoring have been designed and analyzed (either rigorously or heuristically) using our current {\sl knowledge and understanding} of smooth numbers. The same also applies to many primality tests and algorithms for solving the discrete logarithm problem.

Results about the arithmetic structure of ``typical'' integers are therefore of high cryptographic relevance. As most results are already well publicized in the community, we illustrate them by one example (next section) and refer the reader to~\cite{CrPom} for further information.

Nonetheless, new results and applications keep appearing regularly. The works of Croot, Granville, Pemantle \& Tetali~\cite{CGPT} and of Agrawal, Kayal \& Saxena~\cite{AgKaSa} are typical examples.

In~\cite{CGPT} various results about the arithmetic structure of integers
are used to give a very precise analysis of Dixon's factoring algorithm. In~\cite{AgKaSa} results about shifted primes with a large divisor, see Section~\ref{sec:smooth prime}, form the core of the algorithm.

\subsection{Index Calculus in $\F_p^*$}

We start by highlighting the role of smooth numbers in algorithms solving the Discrete Logarithm Problem:

Namely, given two integers $a$ and $b$ and a prime $p$ we consider the problem of computing $k$ (denoted $k = \Dlog_a\, b $) such that $b \equiv a^k \pmod p$ and $0 \le k \le p-2$.

The algorithm is assembled in two steps. We first use a certain (very strong) assumption and then show how to get rid of it.

\paragraph{Initial Assumption:}
Let us fix some $y$ (to be optimized later) and {\sl assume} that we know the discrete logarithms of {\sl all} primes $p_1, \ldots, p_s$ up to $y$ where $s = \pi(y)$.

Under this assumption we perform the following steps:

\begin{description}

\item{{\sl Step~1:}} Pick a random integer $m$ and compute
$$
c \equiv b a^m \equiv a^{k+m} \pmod p, \qquad 0 \le c < p.
$$
Note that
$$\Dlog_a\, c  \equiv \Dlog_a\, b + \Dlog_a\, a^m
\equiv \Dlog_a\, b + m \pmod {p-1}.
$$

The {\sl cost} of this step is negligible.

\item{{\sl Step~2:}}  Try to factor $c$, assuming that $c$, treated as an integer, is $y$-smooth.

Let
$$c = p_1^{\alpha_1} \ldots p_s^{\alpha_s}$$

For doing so use trial division or the elliptic curve factorization algorithm~\cite{Len}.

Note that
$$\Dlog_a\, c \equiv \alpha_1\Dlog_a\, p_1 + \ldots + \alpha_s\Dlog_a\, p_s\pmod {p-1}.$$
The {\sl cost} of this step is about $y$ operations (less if~\cite{Len} is used).

\item{{\sl Step~3:}} If the previous step succeeds, output
$$
\Dlog_a\, b  \equiv \alpha_1\Dlog_a\, p_1 + \ldots + \alpha_s\Dlog_a\, p_s - m \pmod {p-1},
$$
otherwise repeat the first step.

The {\sl cost} of this step is about $p/\Psi(p,y) = u_p^{(1 + o(1))u_p}$ iterations, where
$$
u_p = \frac{\log p}{\log y}
$$
(under the assumption that $c<p$ is a random).

\end{description}

Thus the {\sl total cost}, ignoring nonessential factors, is about $y u_p^{u_p}$.

Taking $y = \exp\(\sqrt{\log p \log \log p}\)$ we get an algorithm of complexity $$\exp\(2 \sqrt{ \log p \log \log p }\)$$

but\ldots it is premature to celebrate the victory, as we need to get rid of the assumption that the discrete logarithms of all small
primes are available.

\paragraph{Removing the Assumption:}
We apply the same algorithm for each $p_i$, $i=1, \ldots, s$ as $b$. Then at Step~3 we get a congruence
$$
\Dlog_a\, p_i \equiv \alpha_{1,i}\Dlog_a\, p_1 + \ldots + \alpha_{s,i}\Dlog_a\, p_s - m_i \pmod {p-1},
$$
for $i=1, \ldots, s$.

We cannot find $\Dlog_a\, p_i$ immediately but after getting such relations for every $p_i$, $i=1, \ldots, s$, we have a system of $s$ linear congruences in $s$ variables. If the system is not of full rank we continue to generate a few relations until a full rank system is reached (this overhead is negligible as most ``random'' matrices are non-singular). Therefore the cost of creating such a system of congruences is about $y^2u_p^{u_p}$ and the cost of solving it is about $y^3$ (lesser if fast linear algebra algorithms are used,  see, for example,~\cite{vzGG}).

Choosing $y$ optimally, we obtain an algorithm of complexity $$\exp\( \cO\(\sqrt{ \log p \log \log p }\)\).$$

The above approach can be improved and optimized in many ways finally yielding a {\sl subexponential} algorithm of asymptotic complexity $$\exp\(\sqrt{(2 + o(1)) \log p \log \log p }\)$$ that can also be rigorously analyzed; this is done by Pomerance in~\cite{Pom}.

We have presented the above example because of its illustrative value although a much faster algorithm exists: the {\sl number field sieve}, see~\cite{CrPom}, whose complexity is $\exp\( \cO\((\log p)^{1/3}(\log \log p )^{2/3}\)\)$.

One can note that the above approach uses both the structure of finite fields and the properties of smooth numbers. Thus a prime field is substantial. Over an extension of a field of small characteristic, such as $\F_{2^n}$, elements can be represented by polynomials and thus smooth polynomials play the role of smooth integers. Hence, the results of Section~\ref{sec:smooth poly} become of great importance.

We note that despite a very common belief that the discrete logarithm problem is solvable in subexponential time, this is {\sl not proved} as we write these lines.

In other words, although over the last decade fast heuristic algorithms for the discrete logarithm problem  have been  designed to work over any finite field, rigorous subexponential algorithms are known only for very specific fields (such as prime fields $\F_p$, their quadratic fields $\F_{p^2}$ or fields $\F_{p^m}$ with a fixed $p$), see~\cite[Section~6.4]{CrPom} for more details.

It is also clear that the above approach does not apply to the discrete logarithm problem in the elliptic curve settings where smoothness admits no analogous notion.

\subsection{Textbook ElGamal Encryption}

The ElGamal cryptosystem~\cite{EG} makes use of two primes $p,q$ with $q \mid p-1$ and an element $g \in \F_p$ of order $q$ (all of which are public),
see also~\cite[Section~8.6]{Buchm},
or~\cite[Sections~8.4.1 and~8.4.2]{MOV},
or~\cite[Section~6.1]{Sti} for further details.

The receiver chooses a random {\sl private key} element $x \in \Z_q$ and computes the {\sl public key} $X = g^x \in \Z_q$.

\begin{description}
\item{{\sl Encryption:}}  To encrypt a message $\mu \in \F_p$, the sender chooses a random $r \in \Z_q$, computes
$R = \mu X^r\in \F_p$, and $Q = g^r\in \F_p$ and sends the pair $\(R,Q\) = \(\mu X^r, g^r\)$.

\item{{\sl Decryption:}} The receiver computes (in $\F_p$)
$$
S = Q^x = g^{xr} = X^r \mand \frac{R}{S} =\frac{R}{X^r} = \mu.
$$
\end{description}

As most public key cryptosystems, the ElGamal protocol is quite slow. It is hence traditionally used to wrap a block-cipher key used for securing the subsequent communication flow.

Doing this in a ``textbook fashion", means that $\mu$ is a rather small integer. For example, $p$ can be about 500 bits long to thwart discrete logarithm calculation attempts, but $\mu$ can be only 80 bits long to resist the brute force search.

Boneh, Joux \& Nguyen~\cite{BoJoNg} have shown that in this case, with a reasonable probability, $\mu$ can be recovered significantly
faster then by any of the above two attacks.

Let  $\cG_q$ be the subgroup of $\F_p^*$ of order $q$ generated by $g$. We note that $R = \mu U$ where $U \in \cG_q$.

Let us assume that $1 \le \mu \le M$ (where $M$ is much smaller than $p$). We also choose some bound  $K$ which is a parameter of the algorithm (controlling the trade-off between complexity and success probability).

\begin{description}

\item{{\sl Step~1:}}  Compute $R^q = \mu^q U^q = \mu^q$.

\item{{\sl Step~2:}}  For $k = 1, \ldots, \rf{K}$ compute, sort and store $k^q$ in a table.

\item{{\sl Step~3:}}  For $m= 1, \ldots, \rf{M/K}$ compute
$$\frac{R^q}{m^q} = \(\frac{\mu}{m}\)^q$$
and check whether this value is present in the table of Step~2.

\item{{\sl Step~4:}}  Output $\mu = km$ if there is a {\sl  match}.
\end{description}

This algorithm always works with $K = M$ (for example, $m = \mu$, $k=1$, which is essentially a form of brute force search).

A better choice is $K = M^{1/2 + \varepsilon}$. Using~\eqref{eq:large int}, we see that the algorithm succeeds for a positive proportion of messages.
That is, it works  because with a sufficiently high probability a random positive integer $\mu \le M$ has a representation $\mu = km$ with $1 \le k, m \le M^{1/2 + \varepsilon}$.

In other words, taking $ M = 2^{80}$ (as in the above example as a standard key size for a private key cryptosystem) we see that the attack runs in a little more than $2^{40}$ steps.

\subsection{Affine RSA Padding}
\label{sec:Fix-Pad attack}

The RSA signature scheme~\cite{RSA} makes use of the following parameters: a composite modulus $N$, a public exponent $e$ and private  exponent $d$
which satisfy the congruence:
$$
ed \equiv 1 \pmod {\varphi(N)},
$$
see also~\cite[Section~8.3]{Buchm},
or~\cite[Section~8.2]{MOV},
or~\cite[Section~5.3]{Sti}.

The {\sl signature} $s\in \Z_N$ of a  message $m \in \Z_N$ is computed as follows: $s \equiv m^d \pmod N$. {\sl Verification} consist in checking
that $m \equiv s^e \pmod N$.

If this is applied in this ``textbook'' form, the scheme becomes susceptible to a {\sl chosen message attack} which works as follows.

Assume that the attacker, wishing to sign a target message $m$, has the ability to ask the legitimate signer to sign seemingly meaningless messages. Then the attacker can:

\begin{itemize}
\item choose a random $m_1$ and compute $$m_2\equiv \frac{m}{m_1} \pmod N;$$

\item query the signatures $s_i \equiv m_i^d \pmod N$ for $i=1,2$ from the legitimate signer;

\item and compute $s \equiv  s_1s_2 \pmod N$.
\end{itemize}

This works because
$$s \equiv s_1s_2 \equiv  m_1^dm_2^d \equiv  (m_1m_2)^d  \equiv  m^d \pmod N.
$$

In other words, because RSA is {\sl homomorphic} with respect to multiplication, a multiplicative relation between messages shadows a similar relation between the signatures.

A natural defense against this attack is to restrict the signature and the verification algorithms to messages of a prescribed {\sl structure}. For example, if $N$ is $n$ bits long, it is requested that the meaningful message part $m$ is only $\ell$ bits long to which a fixed $(n-\ell)$-bit string (called {\sl padding pattern}) is appended. Clearly in the above example $m_1$ can be chosen to comply with this format but $m_2$ is unlikely to fulfill this constraint,
which thwarts the attack.

In the case of affine padding, signed messages have the following structure:

\medskip

\centerline{\fbox{ fixed $(n-\ell)$-bit padding $P$  $|$
  \phantom{M} $\ell$-bit message $m$ \phantom{M}
 }}

\medskip

Thus, denoting $R(m) = P+m$ we see that the signature $s(m)$ of an $\ell$ bit message $m$ is computed as
$$
s(m) \equiv R(m)^d \pmod N, \qquad 1 \le s(m) \le N,
$$
(that is, $P = 2^\ell \Pi$ where $\Pi$ is the appended padding pattern).

In a thread of works by Misarsky~\cite{Mi}, Girault  \&  Misarsky~\cite{GiMi1,GiMi2} and Brier, Clavier, Coron \& Naccache~\cite{BCCN},
{\sl existential forgery} attacks on affine-padded RSA signatures have been  progressively developed and refined.

Lenstra \& Shparlinski~\cite{LenShp} have improved~\cite{BCCN} by redesigning it as a {\sl selective forgery} attack, where the attacker can sign {\sl any} message.

Let us start by presenting the basic technique introduced in~\cite{BCCN}.

Our goal is to find four distinct $\ell$-bit messages  $m_1,m_2,m_3,m_4$  such that
\begin{equation}
\label{eq:R prod}
R(m_1)\cdot R(m_2)\equiv R(m_3)\cdot R(m_4)\pmod N.
\end{equation}
In this case we obtain
$$
s(m_1)\cdot s(m_2)\equiv s(m_3)\cdot s(m_4)\pmod N,
$$
and hence a signature on $m_3$ can be computed from signatures on $m_1, m_2, m_4$.
In~\cite{BCCN} this has been applied to the case where all four messages are considered as variable  $m_1,m_2,m_3,m_4$ (which leads to an existential signature forgery), while in~\cite{LenShp} the message $m_4$ is assumed to be fixed (which leads to a selective signature forgery).

One verifies that the congruence~\eqref{eq:R prod} is equivalent to
$$ P(m_3+ m_4- m_1 - m_2) \equiv m_1m_2 - m_3m_4 \pmod N.$$

With
$$
x=m_1-m_4,\quad y=m_2-m_4,\quad z=m_3+m_4-m_1-m_2
$$
this becomes
\begin{equation}
\label{eq:xyz prod}
(P+m_4)z \equiv xy\pmod N.
\end{equation}
We note that if $m_4$ is already chosen, the values of $x$, $y$ and $z$ define $m_1$, $m_2$ and $m_3$ uniquely.

The congruence~\eqref{eq:xyz prod} is trivial to solve without any restrictions on the variables, but in fact we need ``small'' $x$, $y$ and $z$ about $\ell$ bits long, which is a much harder constraint to deal with.

We show how to solve it when
$$
\ell = \(\frac13 + \varepsilon\) n.
$$

Before we proceed with the algorithm we note  that this choice of $\ell$ is close to the limit of this approach given that for any fixed $\varepsilon>0$ a ``typical'' polynomial congruence in three variables
\begin{equation}
\label{eq:Gen Cong}
F(x,y,z) \equiv  0 \pmod N
\end{equation}
is unlikely to have a integer solution $(x,y,z)$ with
\begin{equation}
\label{eq:small xyz}
1 \le x,y,z \le N^{1/3 - \varepsilon}.
\end{equation}

This is because $F(x,y,z)$ takes only $N^{1-3\varepsilon}$ possible values for such  $x$, $y$ and $z$, thus~\eqref{eq:Gen Cong}
is solvable under the condition~\eqref{eq:small xyz} only with exponentially small ``probability'' of order  $N^{-3\varepsilon}$ (this estimate assumes that $F$ behaves like a random trivariate function and hence must not be taken literally).

Now, to find $\ell$-bit solutions to the congruence~\eqref{eq:xyz prod} we first consider the congruence
\begin{equation}
\label{eq:small sols}
(P+ s)z \equiv w \pmod N,
\end{equation}
where $|s|\le   N^{1/3 + \varepsilon}$ is given and the variables $w$ and $z$ satisfy
$$
w \le N^{2/3 + 2\varepsilon}\mand |z| \le N^{1/3}.
$$

Let ${R_i}/{Q_i}$ denote the $i$-th continued fraction convergent to $(P+s)/{N}$, $i =1, 2, \ldots$.

Then
$$\left| \frac{P+s}{N} - \frac{R_i}{Q_i} \right| \le \frac{1}{Q_iQ_{i+1}}.$$

We now define $j$ by the inequalities $Q_j<N^{1/3}\leq Q_{j+1}$ and set
$$
w=|(P+s)Q_j-NP_j|.
$$
Then
$$0<w\leq \frac{N}{Q_{j+1}}< N^{2/3}\mand
 (P+s)z\equiv w\pmod N$$
for some  $z$ with $|z|<N^{1/3}$, namely  $z = \pm Q_j$.

Now, one can certainly try to  apply the above procedure with $s = m_4$ and then check whether $w$ can be factored as $w = xy$
with $1 \le x,y< 2^\ell$. However, this would unfortunately happen only for a rather sparse sequence of messages $m_4$, so we may
try to randomize the above idea as follows:

\begin{itemize}

\item  Pick a random integer $r$ with $0 \le r < {N^{\varepsilon}}/{2}$ and find $$w \equiv \(P+m_3 - r\fl{N^{1/3}}\)z \pmod N$$ with $1 \le w < N^{2/3}$.

\item Let $u = w + r\fl{N^{1/3}}z$, thus
$$u \equiv (P+m_3 )z \pmod N$$ and $1\le u < N^{2/3+\varepsilon}$ (provided
that $N$ is large enough).

\item Try to factor $u$ using the elliptic curve factoring method which
requires $\exp\(\(2 + o(1)\)\sqrt{\log p \log \log p}\) (\log N)^{O(1)}$
bit operations,
where $p = P_2(u)$. Abort this steps if this takes  longer then some selected time bound.
Thus we abort this steps if $P_2(u)$ is large.

\item Try to find $x,y$ with $u = xy$ and  $1 \le x,y < 2^\ell$.

\item If successful, compute  $m_1, m_2, m_3$, otherwise try another pair $(u,z)$.
\end{itemize}

The above works because eventually we hit a reasonably good $u$ of the form $u = P(u) v$ where $p = P(u) \le 2\sqrt{u}$ and $P(v) = P_2(u)$ is small.

For such a $u$ it is easy to find a representation $u = xy$ with integers $x$ and $y$ of the desired size, $1 \le x, y \le 2^\ell$.

The algorithm seems to be very hard to analyze rigorously, but the heuristic analysis given in~\cite{LenShp} predicts the runtime as $L_{N}(1/3, 1)$ which is substantially faster than $$L_N\(1/3,\(\frac{128}{27}\)^{1/3}\)\approx L_N(1/3,1.68),$$ where
$$
L_N(\alpha,\gamma) =
\exp\((\gamma + o(1)) (\log N)^\alpha (\log\log N)^{1-\alpha}\).
$$
for $N\to\infty$.

Furthermore, Lenstra \& Shparlinski~\cite{LenShp} give a 1024-bit affine padding forgery example while direct factorization of moduli of this size is currently beyond reach.

A challenging open question is to find a way to use more signatures and thereby extend the range of $\ell$ which can be attacked. Recent progress on this question, from a rather unexpected direction, can be found in a work by Joux, Naccache \& Thom\'e~\cite{JNT}.

Finally, although unrelated, other recreative applications of {\sl ad-hoc} factoring in cryptanalysis can be found in~\cite{pkcs} and~\cite{CorNac1}.

\subsection{Desmedt-Odlyzko Attack}
\label{sec:DO attack}

The attacks that we have just described work because of the  inability of the
affine padding to eradicate the  homomorphic properties of RSA. However, there are other attacks that apply in theory to any type of message padding.

In~\cite{DesOd}, Desmedt \& Odlyzko describe an existential RSA signature forgery scenario. Here, the opponent is allowed to query from the legitimate signer $e$-th roots (signatures) of validly padded messages of his choosing. Having done so, the opponent crafts a new validly padded signature of his own on a message left unsigned by the legitimate signer. The attack works as follows:

\begin{description}

\item{{\sl Step~1:}}  Select a bound $y$ and let $p_1, \ldots , p_s$ be the primes  up to $y$, that is, $s =\pi(y)$.

\item{{\sl Step~2:}}  Find $k + 1$ messages $m_i$ which are  $y$-smooth and factor them
$$
m_i = \prod_{j=1}^k p_{j}^{\alpha_{i,j}}
$$
(for example, using elliptic curve factorization \cite{Len} or trial division).

\item{{\sl Step~3:}}  Solve in $u_1, \ldots, u_k \in \{0, \ldots, e-1\}$
$$
\sum_{i=1}^{k}  \alpha_{i,j} u_i \equiv \alpha_{k+1} \pmod e,
\quad j=1, \ldots, k,
$$
and write
$$
\sum_{i=1}^{k}  \alpha_{i,j} u_i + \gamma_i e = \alpha_{k+1} , \quad j=1, \ldots, k.
$$
Thus
$$
m_{k+1} \equiv  r^e \prod_{i=1}^k m_i^{u_i}  \pmod N,
$$
where
$$
r =\prod_{j=1}^k p_{j}^{\gamma_{j}}.
$$

\item{{\sl Step~4:}} Obtain from the legitimate signer the signatures $s_i$ on  $m_i$ for $i  = 1, \ldots, k$ and forge the signature on $m_{k+1}$
as
$$
s  \equiv  r \prod_{i=1}^{k} s_i^{u_i} \equiv r  \prod_{i=1}^k m_i^{d u_i} \pmod N.
$$
\end{description}

We see that $s$ is a {\sl  valid signature} since:
$$
s^e  \equiv   r^e \prod_{i=1}^k m_i^{ed u_i}  \equiv  r^{e} \prod_{i=1}^k m_i^{u_i} \equiv  m_{k+1}^{d}   \pmod N.
$$

Clearly, in a real-life scenario, instead of generating the $y$-smooth
messages $m_1, \ldots, m_k$ on which valid signatures are required, the attacker may also passively monitor the legitimate date exchange and test each signed message for smoothness, putting aside smooth messages together with their signatures.

Using such a strategy, the expected occurrence of a smooth $m$ is about once in
$$\frac{N}{\Psi(N,y)} = u_N^{(1 + o(1))u_N}$$
signature rounds, where
$$
u_N = \frac{\log N}{\log y}
$$
(under the assumption that $m<N$ is random).

Thus, the effort of collecting $y$ messages is $y u_N^{(1 + o(1))u_N}$ which for
$$
y =  \exp\(  \sqrt{0.5 \log N \log \log N }\)
$$
optimizes as $\exp\(\sqrt{(2+ o(1))  \log N \log \log N }\)$.

Coppersmith, Coron, Grieu, Halevi, Jutla, Naccache \& Stern~\cite{CCGHJNS} added a number of improvements and generalizations to this attack and applied it successfully to a number of industry standards.

\subsection{Small Prime Based Public-Key Encryption}

Products of small primes can also be used for public-key encryption.
The idea, due to Naccache \& Stern~\cite{NS1}, is based on the following problem:

\begin{quotation}
{\sl Given a prime $p$, a positive integer  $f<p$ and a set of integers $\{v_1, \ldots, v_n\}$, find a binary vector $x$ such that
$$
f \equiv \prod_{i=1}^n v_i^{x_i} \pmod p,
$$
if such a vector exists.}
\end{quotation}

It is easy to observe that if the $v_1, \ldots, v_n$  are relatively prime and much smaller than $p$, then the exponent vector $x$ can be found in polynomial time by factoring $f$.
Indeed, instances where
$$
p > \prod_{i=1}^n v_i \mand \gcd\(v_i,v_j\)=1, 1 \le i< j\le n,
$$
are easy.

Such an easy instance can be hidden by extracting the $s$-th modular root of each $v_i$, where $s$ is a secret integer with $\gcd(s,p-1)=1$.

More formally, let $p$ be a large public prime and denote by $n$ the
largest integer such that:
$$p>\prod_{i=1}^{n}p_{i}
$$
where $p_i$ is the $i$-th prime.

The secret-key $s < p-1$ is a random integer such that $\gcd(p-1,s)=1$ and the public-keys are
the $s$-th roots:
$$u_i\equiv p_i^r \equiv p_i^{1/s} \pmod p, \ 0 \le u_i < p, \qquad i =1, \ldots, n,
$$
where $r$ satisfies
$$
rs \equiv 1 \pmod {p-1}.
$$

An $n$-bit  message $(m_1, \ldots, m_n)$ is encrypted as
$$
c \equiv \prod_{i=1}^{n}u_i^{m_i}  \pmod p
$$
and recovered by computing
$$
f \equiv c^s \pmod p, \qquad  0 \le t < p,
$$
and then
$$
m_i = \left\{
\begin{array}{ll}
0, & \text{ if } p_i\nmid f, \\
1, & \text{ if } p_i\mid f,
\end{array}
\right.  \qquad i =1, \ldots, n.
$$

We refer the reader to~\cite{NS1} and~\cite{S-MNS} for more information on this somewhat unusual public-key encryption scheme, whose encoding
idea dates back to 1931, see Section~\ref{sec:Godel}.

\subsection{G\"{o}del Numbers}
\label{sec:Godel}

In his famous work published in 1931, G\"{o}del~\cite{Go} uses a mapping of mathematical expressions into integers based on divisibility by small prime factors.

G\"{o}del~\cite{Go} starts by assigning a unique natural number $\tau(\xi)$ to each basic mathematical symbol $\xi$ in the formal language of arithmetic he is dealing with\footnote{for example, $\xi\in\{\exists,\forall,\Rightarrow,+,-,\times,\div,0,1,2,\ldots\}$} (in other words $\tau$ is a symbol-to-integer dictionary).

To encode an entire mathematical {\sl expression} $\Xi$,  which is nothing but an ordered sequence of mathematical symbols:
$$
\Xi= \langle \xi_1,\ldots,\xi_n\rangle.
$$
G\"{o}del~\cite{Go} uses the following system: Each atomic symbol being associated to a positive integer {\sl via} $\tau$, the mathematical expression is mapped into $\N$ as:
\begin{equation}
\label{eq:Godel Numb}
\tau(\Xi) = \prod_{i=1}^{n} p_i^{\tau(\xi_i)}\in\N,
\end{equation}
where $p_i$ stands for the $i$-th prime.

Given that any number obtained this way can be uniquely factored into prime factors, it is possible to effectively and unambiguously recover any mathematical expression $\Xi$ from its G\"{o}del number $\tau(\Xi)$~\eqref{eq:Godel Numb}.

G\"{o}del~\cite{Go} uses this scheme at two levels: first, to encode sequences of symbols representing formulae, and second, to encode sequences of formulae representing proofs.
This has allowed him to show a correspondence between statements about natural numbers and statements about the provability of theorems dealing with natural numbers,  which is the cornerstone of the celebrated G\"{o}del  Incompleteness Theorem~\cite{Go}.

\subsection{Error Correction with Products of Small Primes}

Interestingly, G\"{o}del's encoding~\eqref{eq:Godel Numb} can also be used for error correction.

Error-correcting codes are used to protect information sent over noisy channels against transmission errors. In~\cite{CorNac2,Nac}, Coron \& Naccache describe an unusual error-correcting code based on modular arithmetic.

Let $m$ be the $n$-bit message to encode; we denote by $m_i$ the $i$-th bit of $m$. We let $p_i$ be the $i$-th prime, starting with $p_1=2$. Let $t$ be the number of errors which can be corrected. We generate a prime $p$ such that:
\begin{equation}
\label{equp}
2   p_n^{2   t} \leq p < 4   p_n^{2  t}
\end{equation}
(which, of course, always exists).

Given $m$, we generate the following ``redundancy'':
\begin{equation}
\label{defcm}
c(m) \equiv \prod_{i=1}^n p_i^{m_i} \pmod p, \qquad 0 \le c(m) < p.
\end{equation}

The integer $c(m)$ is protected by using an error-correcting code $\mu$ resilient to $t$ transmission errors.

The encoded message $E(m)$ is defined as $E(m)=\langle m,\mu(c(m))\rangle$.

Let $\tilde E(m)$ be the received version of $E(m)$ where at most $t$ errors occurred:
$$
\tilde E(m)  =   E(m) \xor e
$$
where $e$ is an error vector of Hamming weight at most $t$, and $\xor$ stands for bit-wise addition.

Splitting $e = \langle e_m,e_c\rangle $ into parts corresponding to errors in $m$ and to errors in $\mu(c(m))$, we obtain
$$
\tilde E(m) =  \langle \tilde{m}, \tilde c(m) \rangle
 = \langle m \xor e_m , \mu(c(m)) \xor e_c \rangle.
$$

Since $\mu$ can correct $t$ errors, and $e_c$ is only a part of $e$ (whose {\sl total} Hamming weight is $t$), $c(m)$ can be safely recovered from $\mu(c(m)) \xor e_c$.

The receiver computes:
$$
s\equiv \frac{c(\tilde{m})}{c(m)} \equiv \frac{c(m \xor e_m)}{c(m)} \pmod p.
$$

Using~\eqref{defcm} the integer $s$ can be written as:
$$s \equiv a/b \pmod p,
$$
where
$$
a =  \prod_{\substack{i =1,\ldots,n\\ \widetilde m_i=1\\ m_i=0}} p_i  \mand
b =  \prod_{\substack{i =1,\ldots,n\\ \widetilde m_i=0\\ m_i=1}} p_i .
$$
Since $\tilde{m}$ suffered at most $t$ errors, we have
$$
\max\{a,b\} < p_n^t .
$$

A result of Stern, Fouque \& Wackers~\cite{SFW} shows that given $s$ one can recover $a$ and $b$ efficiently.
The algorithm is based on the Gauss reduction algorithm for finding the shortest vector in a two-dimensional lattice~\cite{vallee}.  More precisely, let $p$  be an prime with  $p > 2AB$ for some $A\in\R$ and $B\in\R$. Let $a,b \in\Z$ be such that $|a| \le A$ and $0<b \le B$. Then given $p$, $A$, $B$ and $s \equiv a b^{-1} \pmod p$, one can recover $a$ and $b$ in polynomial time. Note that the condition $p > 2AB$ guarantees the uniqueness of $a$ and $b$. A very similar argument has been used in Section~\ref{sec:Fix-Pad attack} to find small solutions to~\eqref{eq:small sols}.

Taking $A=B=p_n^t-1$, we have from~\eqref{equp} that $2AB<p$.
Moreover, $0 \leq a \leq A$ and $0 <b \leq B$. Therefore, we can recover $a$ and $b$ from $s$ in polynomial time. By testing the divisibility of $a$ and $b$ by the small primes $p_1, \ldots, p_n$, one can recover $e_m = \widetilde m \xor m$ and hence $m = \widetilde m \xor e_m$.

The process assumes the existence of an error correcting code $\mu$. Note that $\mu$ can be nothing but the procedure that we have just presented in miniature. In other words, the described encoding procedure can be iterated to protect $c(m)$ using a new, much smaller, set of primes. In turn, yet another encoding iteration is used at the third level of encoding and so on. Finally, the smallest and last layer can be protected by simple replication ($2t+1$ times) and decoded using a majority vote.

The proposed code turns out to provide efficient decoding for some specific parameter combinations. For instance, denoting by $\mu$ Reed-Muller encoding, and assuming that 5812-bit messages need to be protected against 31 transmission errors, the size of $\mu(m)$ is 8192 bits, whereas the hybrid encoding $\langle  m,\mu(c(m))\rangle$ is only 7860 bits long.

More examples and details can be found in~\cite{CorNac2,Nac}. The full asymptotic analysis of this scheme still remains to be worked out.

\subsection{Private Information Retrieval with Products of\break Small Primes}

A {\sl Private Information Retrieval} (PIR) scheme is a combination of encoding and encryption which allows a user to retrieve the $k$-th bit of an $n$-bit database, without revealing to the database owner the value of $k$.

Gentry \& Ramzan~\cite{GeRa} have used the Chinese Remainder Theorem and properties of products of small primes to design a PIR scheme. The construction of~\cite{GeRa} requires a cyclic group $\cG$ whose order $t = \# \cG$ has a prescribed arithmetic structure; namely a product of a large prime and a very smooth integers. This makes the results of~\cite{BaSh1,PomShp,Ten4} relevant to this problem, see also Section~\ref{sec:Smooth Div}.

\subsection{Zero-Knowledge with Products of Small Primes}

A zero-knowledge proof (ZKP) is a protocol allowing Alice to convince Bob that she knows a secret $s$ without revealing to Bob information on $s$.

The best-known ZKP is probably the protocol of Fiat \& Shamir~\cite{FS} which uses an RSA modulus $N$ and $k$ quadratic residues $v_i$ as public parameters. In its simplest version, Alice uses the $k$ modular square roots $s_i$ such that $s_i^2 \equiv v_i \pmod N$ as secret identification keys. The protocol is:

\begin{itemize}
\item Alice picks a random $r\in\Z_N$ and sends to Bob $x \equiv r^2 \pmod N$.
\item Bob picks a random binary vector
$$
e=\langle e_0,\ldots,e_{k-1}\rangle
$$
and sends it to Alice.
\item Alice replies to Bob with:
$$
y \equiv r\prod_{i=0}^{k-1} s_i^{e_i} \pmod N.
$$
\item  Bob verifies that:
$$
y^2 \equiv x\prod_{i=0}^{k-1} v_i^{e_i} \pmod N.
$$
\end{itemize}

To ease Bob's computational burden, Micali \& Shamir~\cite{MS} suggest to use very small $v_i$-values. As it turns out, using small primes as $v_i$-values presents particular security and simplicity advantages.

\subsection{The Generalized Diffie-Hellman Problem}

Recently, several cryptographic schemes based on the following assumption appeared:

Let $g$ be an element of prime order $p$ of a ``generic'' Abelian group $\cG$. That is, we assume that $\cG$ is a group where only ``generic'' attacks, such as Shanks' or Pollard's algorithms exist and take about $\sqrt{p}$ operations,
see~\cite[Sections~10.3 and~10.4]{Buchm},
or~~\cite[Sections~3.6.2 and~3.6.3]{MOV},
or~\cite[Sections~6.2.1 and~6.2.2]{Sti}. For example, one may regard $\cG$ as the group of points on an elliptic curve over a finite field.

%%igor Period, Solving
The {\sl traditional Diffie-Hellman problem} is defined as follows:
\begin{quotation}
{\sl Given $g^x$ and $g^y$, compute $g^{xy}$.}
\end{quotation}
Solving this problem is believed to be hard.

Due to the identity
$$
g^{2xy} = g^{(x+y)^2} g^{-x^2} g^{-y^2}
$$
and the fact that computing square roots in groups of prime order is easy the Diffie-Hellman problem can be reformulated in a shorter form:
\begin{quotation}
{\sl Given $g^x$, compute $g^{x^2}$.}
\end{quotation}

On the other hand, many cryptographic protocols rely on the presumed hardness of the following {\sl generalized Diffie-Hellman problem}:
\begin{quotation}
{\sl Given $n$ powers $g^x, \ldots g^{x^{n}}$, compute $g^{x^{n+1}}$.}
\end{quotation}

Intuitively it may seem that, despite the fact that more information on $x$ leaks out in the generalized Diffie-Hellman 
%%igor setting, 
settings solving it is not easier than solving the  
traditional Diffie-Hellman problem with the same parameters.

Surprisingly, Brown \& Gallant~\cite{BrGal} and
Cheon~\cite{Cheon}, have shown this intuition to be wrong.

Here are some results of Cheon~\cite{Cheon}:

%%igor . -> ;
\begin{itemize}
\item given  $g^x$ and $g^{x^d}$ for some $d \mid p-1$, one can find $x$ in time about $\cO\(\sqrt{p/d} + \sqrt{d}\)$ (which is $\cO\(p^{1/4}\)$ for $d \sim \sqrt{p}$);

\item  given  $g^x, \ldots g^{x^d}$ for some $d \mid p+1$, one can find $x$ in time about $\cO\(\sqrt{p/d} + d\)$
(which is $\cO\(p^{1/3}\)$ for $d \sim p^{1/3}$).
\end{itemize}

This brings up the question of estimating the probability at which primes $p$ are such that $p\pm 1$ has a divisor $d$ of a given size.

More specifically,  how rare are primes $p$ such that $p\pm 1$ has a divisor $d \in [n^{1- \varepsilon}, n]$? (which  guarantees the asymptotically
best advantage if we are given $g^x, \ldots g^{x^{n}}$ with $n$ which is not too large).

By the result of Ford from Section~\ref{sec: div shift p} we know that for every $\varepsilon > 0$ this happens for a positive proportion of primes $p$.

Therefore, we conclude that the attack of~\cite{Cheon} can be applied in its asymptotically strongest form with a positive probability.
In other words, the generalized Diffie-Hellman problem {\sl is easier} than the traditional Diffie-Hellman problem.

%%igor
In practical scenarios probably only small values of $d$ can be used.
In this case the bound~\eqref{eq:Div p-1-Brun} can be applied.

\subsection{Large Subgroup Attack}
\label{sec:HMQV}

The Digital Signature Algorithm uses two large primes
$p$ and $q$ such that $q\mid (p-1)$, see~\cite[Section~12.6]{Buchm},
or~\cite[Section~11.5.1]{MOV},
or~\cite[Section~7.4.2]{Sti},

Suppose that $p$ and $q$ are generated using the following straightforward method:

\begin{itemize}
\item select a random $m$-bit prime $q$;

\item randomly generate $k$-bit integers $n$ until a prime
$p=2nq+1$ is reached.
\end{itemize}

%%igor
In~\cite{Men}, Menezes introduces the {\sl Large Subgroup Attack}
on some cryptographic protocols, including a version of the HMQV
protocol, see also~\cite{MenUst}.

The attack can be applied if $n =({p-1})/({2q})$ has a smooth divisor $s > q$. Some upper bounds on the density of such primes with a large smooth divisor are given by
Pomerance \& Shparlinski~\cite{PomShp}.

However this result does not take into account the special structure of $p$ (for example, the presence of a large prime divisor $q \mid (p-1)$), so it does not (quite) apply.

Furthermore, in the above situation lower bounds become more important. Determining such bounds is unfortunately a much harder question.

On the other hand, using the results of Banks \& Shparlinski~\cite{BaSh1} and Tenenbaum~\cite{Ten4} mentioned in Section~\ref{sec:Smooth Div}
one can get an estimate of the probability $\eta(k,\ell, m)$ that a $k$-bit integer $n$ has a divisor $s > 2^m$ which is $2^\ell$-smooth. Then, assuming that shifted primes $p-1$ behave like ``random'' integers, one can address the original question (at least heuristically).

The most interesting choice of parameters as we write these lines is:
$$
k = 863, \qquad m = 160, \qquad \ell = 80
$$
(which produces a $1024$-bit prime $p$).

It has been shown in~\cite{BaSh1} that for these parameters, the theoretic estimates (together with some heuristic assumptions about the distribution
of primes in the sequence $2qn+1$ for $n$ having a large smooth part) suggest  that the attack succeeds  with probability
$$
\eta(863,80,160)\approx 0.09576> 9.5\%
$$
over the choices of $p$ and $q$.

We also note that similar attacks on the ElGamal signature scheme and the Diffie Hellman key exchange protocols, have been outlined by  Anderson \& Vaudenay~\cite{AnVa}.

\subsection{Smooth Orders}

Let $l(n)$ be the multiplicative order of 2 modulo $n$, $\gcd(2,n)=1$ (in the following 2 can be replaces by any integer $a \ne 0, \pm 1$).

Motivated by several cryptographic applications, Pomerance \& Shparlinski~\cite{PomShp} has studied the smoothness of $l(n)$ on integers and on shifted primes $n = p-1$. This arises from the desire to clarify whether $g = 2$ can safely serve as an exponentiation base in
discrete logarithm based cryptosystems\footnote{Small values of $g$ allow to significantly speed up square-and-multiply exponentiation.}. However, in order to avoid the {\sl Pohlig-Hellman attack} $l(n)$ must not be smooth, see~\cite[Section~10.5]{Buchm},
or~\cite[Section~3.6.4]{MOV},
or~\cite[Section~6.2.3]{Sti}.

Also, Boneh \& Venkatesan~\cite{BonVen} have shown that the  Diffie-Hellman protocol with the $g=2$ has some additional attractive  bit security properties which are not known for other $g$ values.

Finally we recall that Pollard's $(p-1)$-factorization method works better
when $p \mid n$ features a smooth $l(p)$,
see~\cite[Section~9.2]{Buchm}, or \cite[Section~5.4]{CrPom},
 or~\cite[Section~3.2.3]{MOV}, or~\cite[Section~5.6.1]{Sti} for details.
Some improvements of this algorithm have recently been suggested by Zralek~\cite{Zra1}.

Let us define the following counting functions:
$$L(x,y) = \# \{ p \le x \ : \  l(p) \text{ is $y$-smooth}\}
$$
and
$$N(x,y) = \# \{ n \le x \ : \  l(n) \text{ is $y$-smooth}\}.$$

Pomerance \& Shparlinski~\cite{PomShp} have shown that for
$$\exp\(\sqrt{\log x\log\log x}\,\)\le y \le x,$$
we have
$$L(x,y)\ll u\;\rho\left(\frac{u}{2}\right)\pi(x).$$

It is also noticed in~\cite{PomShp} that it seems quite plausible that in fact the bound also holds with $\rho(u)$ instead of $\rho(u/2)$,
which means that the values of $l(p)$ behave as ``random'' integers.

In fact, this may even happen to be provable under the Generalized Riemann Hypothesis. However this has not been worked out yet and remains an interesting open question.

Furthermore, Banks, Friedlander, Pomerance \& Shparlinski~\cite{BFPS} proved that for
$$\exp\(\sqrt{\log x\log\log x}\,\)\le y\le x$$
we have
$$N(x,y)\le x\exp\(-\(\frac12 + o(1)\)\,u \log \log u\).$$

As in the case of $L(x,y)$, one may expect that the same bound should hold with $1$ instead of $1/2$ in the exponent, but the appearance of
$\log \log u$ instead of $\log u$ seems to be right, see also Section~\ref{sec:smooth phi}.

\subsection{Smooth-Order Based Public Key Encryption}

Smooth orders can also be used constructively to provide public key encryption. Here is one such suggestion due to Naccache \& Stern~\cite{NS2}:

\paragraph{Parameter Generation:}
Let $s$ be a odd, squarefree, $y$-smooth integer, where $y$ is a certain small parameter and let $N=pq$ be an RSA modulus such that
$$
s \mid \varphi(n) \mand \gcd\(s,\frac{\varphi(N)}{s}\)=1.
$$

Typically, we think of $y$ as being a 10 bit integer and consider $N$ to be at least $768$ bits long. Let $g$ be an element whose multiplicative order modulo $N$ is a large multiple of $s$. Publish $N$, $g$ and keep $p$, $q$ and  $s$ secret (note that there are very few possibilities for $s$ so its revealing does not give any dramatic advantage to the attacker).

Generation of the modulus appears rather straightforward: pick a family $p_1 < \ldots < p_k$ of $2k$ small odd distinct primes and set:
\begin{equation}
\label{eq:uvs}
u = \prod_{i=1}^{k} p_{2i-1}, \qquad v =  \prod_{i=1}^{k} p_{2i-1}, \qquad
s = uv = \prod_{i=1}^{k} p_i
\end{equation}
(thus $s$ is $p_k$-smooth).

Find (using trials and primality testing) two large primes $\ell$ and $r$ such that both $p = 2\ell u +1$ and $q = 2r v +1$ are prime and let $N=pq$.

Note that much faster key generation procedures exist, we refer the reader to~\cite{NS2} for more details.

To generate $g$, one can choose it at random in $\Z_N$ and check whether it has the  possible order $\varphi(N)/4$ or $\varphi(N)/2$  modulo $N$,
Note that for any $N$ multiplicative  orders of elements of $\Z_N$ are divisors of the Carmichael function $\lambda(N)$; in the above case $\lambda(N)=\varphi(N)/2$.

The main point is to ensure that $g$ is not a $p_i$-th power modulo $N$ for each $i= 1, \ldots, k$ by testing that
$$
g^{\varphi(n)/p_i} \not \equiv 1 \pmod N, \qquad i= 1, \ldots, k.
$$
The success probability is:
$$\rho= \prod_{i=1}^{k} \(1- \frac{1}{p_i}\).
$$
If the $p_1, \ldots, p_k$ are the first $k$ odd primes, this in turn can be estimated by the Mertens formula as  $\rho\sim 1/\log k$. Another method consists in choosing, for each index $i \leq k$, a random $g_i$ until it is not a $p_i$-th power. With overwhelming probability
$$g = \prod_{i=1}^{k} g_i^{s/p_i}$$
has the multiplicative order at least $\varphi(n)/4$.

\paragraph{Encryption:}
A message $m < s$ is encrypted as
$$c \equiv g^m \pmod N.$$

\paragraph{Decryption:}
The algorithm computes the value $m_i$ of the residue of $m$ modulo each prime factor $p_i$, $i= 1, \ldots, k$ of $s$ given by~\eqref{eq:uvs}, and recovers the message by the Chinese Remainder Theorem, following an idea of Pohlig-Hellman~\cite{PoHe}, see also~\cite{CrPom,MOV}.

Now for every $i= 1, \ldots, k$, to find $m_i$, given the ciphertext $c\equiv g^m \pmod N$, the algorithm computes
$$c_i \equiv c^{{\varphi(n)}/{p_i}}
\equiv g^{{m \varphi(n)}/{p_i}}\equiv g^{{m_i \varphi(n)}/{p_i}} \pmod N, $$
where the congruence $m \equiv m_i \pmod {p_i}$ is used at the last step.

By comparing this result with all possible powers
$$
g^{{j \varphi(n)}/{p_i}}, \qquad j =0, \ldots, p_i-1,
$$
the algorithm finds out the correct value of $m_i$.

The basic operation used by this (non-optimized) algorithm is a modular exponentiation of complexity $\cO\((\log N)^3\)$, repeated at most :
$$k  p_k  \ll  k^2 \log k \ll (\log N)^2  \log\log N$$ times.
Decryption therefore takes $\cO\((\log N)^5 \log\log N \)$ bit operations.

We refer the reader to~\cite{NS2} for more details and optimizations.

%%igor
\subsection{Oracle-Assisted Integer Factorization}

Maurer~\cite{Maur} has designed an algorithm which for any $\varepsilon$, given an integer $N$, requests at most $\varepsilon \log N$ bits of information and factors $N$ in polynomial time.

Unfortunately a rigorous analysis of this algorithm requires
very precise results about the distribution of smooth numbers
in short intervals which currently seems to be beyond reach.
Accordingly, the main result of~\cite{Maur} is conditional
and relies on heuristic assumptions.

\subsection{Pratt Trees}

Highly critical security applications sometimes require {\sl primality proofs}. Here is a way to provide such proofs, due to Pratt~\cite{Pratt}

\begin{itemize}
\item  Check that the would-be prime $p$ is not a perfect power. This is easy, see, for example,~\cite{Bern2,BLP}.

\item   Produce a primitive root $g$ modulo $p$ and provide a proof of this. For that sake it is enough to verify that
$$
g^{p-1} \not \equiv 1 \pmod p \mand g^{(p-1)/q} \not \equiv 1 \pmod p
$$
for all prime divisors $q \mid (p-1)$, so the list of these primes $q$ must also be supplied.

\item  Give a proof that each such $q$ is prime by iterating the above procedure.
\end{itemize}

The whole algorithm can be viewed as a tree, called the {\sl Pratt Tree}, where each node contains a prime (with $p$ as a root) and with $2$ at each leaf.

The algorithm runs in polynomial time and in particular  shows that the decision problem PRIMES is in the complexity class NP (which is
not so exciting nowadays given that, thanks to~\cite{AgKaSa}, we know that PRIMES is actually in P).

Pratt~\cite{Pratt} has shown that the number of multiplications required by this algorithm is $\cO((\log p)^2)$.
On the other hand, Bayless~\cite{Bay} shows that this number is at least $C\log p$ for any fixed $C > 1$ and for almost all primes $p$.

There are, however, many other interesting questions about this tree, such as estimating its height, number of nodes, number of leaves, and so on, in extreme cases and also for almost all primes.

For example, it is obvious that the Pratt Tree's height $H(p)$ satisfies the inequality
\begin{equation}
\label{eq:H upper}
H(p) \ll \log p.
\end{equation}
One can also infer from more general results of K{\'a}tai~\cite{Kat} that for some constant $c > 0$ the inequality
\begin{equation}
\label{eq:H lower}
H(p)  \ge c \log \log p
\end{equation}
holds for almost all primes $p$.

Ford, Konyagin \& Luca~\cite{FKL} have recently given a heuristic argument
suggesting that
$$
H(p) \gg \frac{\log p}{\log \log p}
$$
for infinitely many primes $p$ and also a rigorous proof that
$$
H(p) \ll (\log p)^{0.9622}
$$
holds for  almost all $p$.
It is also shown in~\cite{FKL}  that~\eqref{eq:H lower}
holds for almost all primes $p$ with any
$$
c < \frac{1}{1 + \log 2}.
$$
It seems that the lower bound~\eqref{eq:H lower} is of the right order of magnitude and in fact some heuristic arguments, given in~\cite{FKL},
lead to the conjecture that
$$
H(p) = e \log \log p + \cO(\log \log \log p)
$$
for almost all primes $p$.

A number of other challenging open questions and conjectures can be found in~\cite{FKL}.

Studying other characteristics of the Pratt Tree is also an interesting and little-researched open question.

For instance, Banks \& Shparlinski~\cite{BaSh2} have shown that the length $L(p)$ of the chain $p\mapsto P(p-1)$ satisfies
\begin{equation}
\label{eq:L lower}
L(p) \ge (1 + o(1) ) \frac{\log \log p}{\log \log \log p} \end{equation}
for almost all primes $p$. This corresponds to a particular path in the Pratt Tree. Furthermore, it may be natural to expect that this should actually be the longest path for almost all primes, so it is possible that
$$
L(p) =(1+o(1)) H(p)
$$
for almost all primes $p$.
On the other hand, it seems that  $L(p) < H(p)$ holds
for almost all primes $p$.
Clarifying the matter is an important research challenge. As a first step one may for instance try to use the methods of~\cite{FKL,Kat} to improve~\eqref{eq:L lower} up to the level of~\eqref{eq:H lower}.

\subsection{Strong Primes}

A prime $p$ is called {\sl strong} if $p-1$ and $p+1$ have a large prime divisor, and  $p-1$ has a prime divisor $r$ such that $r-1$ has a large prime divisor,
see~\cite[Section~4.4.2]{MOV}.

To make this definition more formal we say that $p$ is {\sl $y$-strong\/} if $p+1$ has a large prime divisor $q \ge y$, and  $p-1$ has a prime divisor $r$ such that $r-1$ has a prime divisor $\ell \ge y$.

We note that the combination of~\cite{BFPS} and~\cite{PomShp} (see Sections~\ref{sec:smooth prime} and~\ref{sec:smooth phi}) implies that almost all primes are $y$-strong as $\log x/\log y \to \infty$.

Indeed, the cardinality of the set of primes $p\le x$ such that $p+1$ is $y$-smooth is exactly the function $\pi_{1}(x,y)$ discussed in Section~\ref{sec:smooth prime}.

From the set of remaining primes $p\le x$ we remove those for which $p-1$ is divisible by $r^2$ for a prime $r\ge y$.
Since the number of primes $p \le x$ with $p \equiv 1 \pmod{r^2}$
is at most $x/r^2$, the cardinality of this set can be estimated trivially as
$$
\sum_{r \ge y}  \frac{x}{r^2} = \cO\(x/y\).
$$

 Hence it is easy  to see that if one of the remaining primes is not $y$-strong then $\varphi(p-1)$ is $y$-smooth and thus the bounds of $\Pi_{-1}(x,y)$ from Section~\ref{sec:smooth phi} can be applied.

\subsection{Small Prime Based Hash Functions}

The {\sl Very Smooth Hash\/} function, VSH, recently introduced and studied
 by Contini,  Lenstra \& Steinfeld~\cite{CLS}, is defined as follows.

Let $p_i$ denote the $i$-th prime number and let
$$
Q_k = \prod_{i=1}^{k} p_i
$$
denote the product of the first $k$ primes.

Assume that integers $k$ and $N$ satisfy
\begin{equation}
\label{eq: k and n}
Q_k   < N \le Q_{k+1}.
\end{equation}

Let the message length $\ell < 2^{k}$ be a positive integer whose $k$-bit representation (including all leading zeros) is $\ell = \lambda_1\ldots\lambda_{k}$ that is
$$
\ell = \sum_{i=1}^k \lambda_i 2^{i-1}.
$$

The VSH takes an $\ell$-bit message $m = \mu_1, \ldots, \mu_\ell$ and hashes it (in a very efficient way, via a simple iterative procedure) to
$$
h_N(m) \equiv   \prod_{i =1}^{k} p_i^{e_i} \pmod N, \qquad
 0 \le h_N(m)  < N,
$$
where $L = \rf{\ell/k}$, $\mu_s = 0$, for $\ell < s \le Lk$, $\mu_{Lk+i} = \lambda_i$, for $1 \le i \le k$, and
%%\begin{equation}
%%\label{eq: ei's}
$$
e_i = \sum_{j=0}^L \mu_{jk + i} 2^{L-j}, \qquad i =1, \ldots, k.
$$
%% \end{equation}

It is demonstrated in~\cite{CLS} that the VSH also admits a rigorous collision-resistance proof based on ``natural" number theoretic problems which are presumably hard. As the above problem is related to factoring, it is natural to choose $N = pq$ to be an RSA modulus. The design and the suggested parameter choice are both based on classical facts about the distribution of smooth numbers.

In~\cite{BlShp}, Blake \& Shparlinski harness results about the distribution of smooth numbers to provide rigorous support in favor of the security and the distribution properties of
the  VSH. In particular, \cite{BlShp} shows that for almost all RSA moduli and any integer $a$, the probability that for a random  $\ell$-bit message $m$ we have $h_N(m) \equiv a \pmod  N$, is negligible for sufficiently large  values of $\ell$.

This bounds the collision probability and also the probability of finding a second pre-image by brute force.

The above and several other results in~\cite{BlShp} are based on the study of the multiplicative subgroup of $\Z_N^*$ generated by $p_1, \ldots, p_{k}$ for integers $N = pq$ where $p$ and $q$ are distinct primes, satisfying the inequality~\eqref{eq: k and n}.

\section{Conclusion}

Our goal has been to position this paper at the crossroads of cryptography and number theory. We hope that while reading it cryptographers have enriched their arsenal with a large gamut of little-used, yet powerful, number-theoretic methods and results extending beyond the classical facts used in cryptology. On the other hand, it is our hope that number theorists
have  enjoyed learning how smooth numbers can be harnessed to provide encryption, private information retrieval, identification, error correction, hashing, primality proofs and other cryptographic functions.
Furthermore, final tuning and adjusting already known results and techniques may lead to new advances of intrinsic mathematical interest. Our outline, somewhat sketchy and simplified and also sometimes ignoring subtleties, cannot replace a careful and systematic reading of the original number theoretic and cryptographic literature,
such as~\cite{CrPom,HalbRich,Harm2,IwKow,Ten1} and~\cite{Buch,MOV,Sti}, respectively.

\section{Acknowledgements}
%
%The author is very grateful to the organisers of the
% 4th China-Japan Number Theory Conference
%Jianya Liu and Shigeru Kanemitsu for their kind invitation
%to this meeting  and help
%with preparation of this manuscript.

The authors would like to thank Kevin Ford
for   many valuable suggestions.

This second author work was supported in part by ARC grant DP0556431.

\end{document}